\documentclass[12pt,reqno]{amsart}
\usepackage{mathrsfs,amsfonts,amsmath}
\usepackage{lineno}
\numberwithin{equation}{section}

\newtheorem{theorem}{Theorem}[section]
\newtheorem{corollary}{Corollary}[section]

\newtheorem{lemma}{Lemma}[section]
\newtheorem{proposition}{Proposition}[section]
\newtheorem{definition}{Definition}[section]

 \allowdisplaybreaks

\begin{document}
\title[\ Convergence of a quantum lattice Boltzmann scheme]
      {Convergence of a quantum lattice Boltzmann scheme to the nonlinear  Dirac  equation for Gross-Neveu model in $1+1$ dimensions}
\author[N. Li, J. Zhang and Y. Zhang]{\scshape Ningning Li, Jing Zhang \and Yongqian Zhang \\  \,  \\ {\scriptsize  School of Mathematical Sciences  \\ Fudan University, Shanghai 200433, China} }

\address[N. Li]{\small School of Mathematical Sciences,
                                   Fudan University, Shanghai 200433, P. R. China}
\email{\tt nnli20@fudan.edu.cn}
\address[J. Zhang]{\small School of Mathematical Sciences,
                                  Fudan University, Shanghai 200433, P. R. China}
\email{\tt zjgjzx521@163.com}
\address[Y. Zhang]{\small School of Mathematical Sciences,
                                   Fudan University, Shanghai 200433, P. R. China}
\email{\tt yongqianz@fudan.edu.cn}
\let\thefootnote\relax\footnotetext{ Corresponding author: Yongqian Zhang, yongqianz@fudan.edu.cn}


\begin{abstract}
This paper studies the quantum lattice Boltzmann scheme for the nonlinear Dirac equations for Gross-Neveu model in $1+1$ dimensions. The initial data for the scheme are assumed to be convergent in $L^2$. Then  for any $T\ge 0$ the corresponding  solutions for the quantum lattice Boltzmann scheme are shown to be convergent in $C([0,T];L^2(R^1))$ to the strong solution to the nonlinear Dirac equations as the mesh sizes converge to zero. In the proof, at first a Glimm type functional is introduce to establish the stablity estimates for the difference between two solutions for the corresponding quantum lattice Boltzmann scheme,  which leads to the compactness of the set of the solutions for the quantum lattice Boltzmann scheme. Finally the limit of any convergent subsequence of the solutions for the quantum lattice Boltzmann scheme is shown to coincide with the strong solution to a Cauchy problem for the nonlinear Dirac equations.
\end{abstract}
\date{\today}
\maketitle

\section{Introduction}\label{section-intoduc}
The nonlinear Dirac equations for Gross-neveu model in $R^{1+1}$ can be written as
\begin{equation}\label{eq-dirac}
\left\{ \begin{array}{l} u_t+u_x=imv+iN_1(u,v), \\ v_t-v_x=imu+iN_2(u,v),
\end{array}
\right.
\end{equation}
where $(t,x)\in R^2$, $(u,v)\in \mathbf{C}^2$. The nonlinear terms  take the following form
\begin{equation}\label{eq-nonlinearstruc}
N_1=\partial_{\overline{u}}W(u,v), \quad N_2=\partial_{\overline{v}}W(u,v)
\end{equation}
with
\[
 W(u,v)=\alpha |u|^2|v|^2+\beta (\overline{u}v+u\overline{v})^2,
 \]
see \cite{pelinovsky}. Here $\alpha,\beta\in R^1$ and $\overline{u}, \overline{v}$ are complex conjugate of $u$ and $v$. The initial data is given as follows,
\begin{equation}\label{eq-dirac-initialv}
(u, v)|_{t=0}=(u_0(x), v_0(x)).
\end{equation}
(\ref{eq-dirac}) is called  Thirring equation for $\alpha=1$ and $\beta=0$, while it is called Gross-Neveu equation for $\alpha=0$ and $\beta=1/4$; see for instance \cite{thirring,gross-neveu,pelinovsky}. Such model arises in the study of quantum mechanics and general relativity (\cite{gross-neveu} and \cite{thirring}). There have been many works devoted to the local and global well-posedness of Cauchy problem for (\ref{eq-dirac}) in different kinds of Sobolev spaces, see for instance, \cite{bournaveas-zouraris,cacciafesta,candy,deldado,escobedo,huh,huh-moon,huh2,pelinovsky,selberg,zhang,zhang-zhao} and the references therein. For the case that initial data $(u_0(x), v_0(x))\in L^2(R^1)$,  it has been proved in \cite{zhang-zhao} that (\ref{eq-dirac}) and (\ref{eq-dirac-initialv}) has a unique global strong solution in $L^2$.

In this paper we are concerned with a difference scheme called the quantum lattice Boltzmann scheme for (\ref{eq-dirac}) and (\ref{eq-dirac-initialv}) with $(u_0,v_0)\in L^2(R^1)$. Such a scheme, denoted briefly by QLB,  is proposed by  S. Palpacelli, P. Romatschke and S. Succi \cite{palpaceli-r-succi} for (\ref{eq-dirac}), see also Succi and Benzi \cite{succi-benzi,succi2001} for the QLB schemes for linear Dirac equations. The corresponding scheme for (\ref{eq-dirac}) in \cite{palpaceli-r-succi} are given as follows,
\begin{equation}\label{eq-LBE-A}
\left\{ \begin{array}{l}\displaystyle \frac{\widehat{u}^{(h)}-u^{(h)}}{h}=im \frac{\widehat{v}^{(h)}+v^{(h)}}{2}+i\alpha \frac{\widehat{u}^{(h)}+u^{(h)}}{2}|v^{(h)}|^2+i\beta\frac{\widehat{v}^{(h)}+v^{(h)}}{2}G^{(h)} \\ \, \\ \displaystyle
\frac{\widehat{v}^{(h)}-v^{(h)}}{h}=im \frac{\widehat{u}^{(h)}+u^{(h)}}{2}+ i\alpha \frac{\widehat{v}^{(h)}+v^{(h)}}{2} |v^{(h)}|^2+ i\beta\frac{\widehat{u}^{(h)}+u^{(h)}}{2}G^{(h)},
\end{array}
\right.
\end{equation}
where
\[ (u^{(h)},v^{(h)})=(u^{(h)}(x,t),v^{(h)}(x,t)),\]
\[ (\widehat{u}^{(h)},\widehat{v}^{(h)})=(u^{(h)}(x+h,t+h),v^{(h)}(x-h,t+h))\] and
\[ G^{(h)}=G(u^{(h)},v^{(h)})(x,t):= u^{(h)}(x,t)\overline{v^{(h)}(x,t)}+\overline{u^{(h)}(x,t)}v^{(h)}(x,t).\]
Here \[G(u,v)=\overline{u}v+u\overline{v},\]
  and the function $(u^{(h)},v^{(h)}) $ is piecewise-constants valued, and  satisfies
\begin{equation}\label{eq-funct-value}
(u^{(h)}(x,t),v^{(h)}(x,t))=(u^k_n,v^k_n), \, \, (x,t)\in [nh,(n+1)h)\times [kh,(k+1)h)
\end{equation}
for the any integers $n$ and $k\ge 0$, where
\[(u^k_n,v^k_n)=(u^{(h)}(nh,kh),v^{(h)}(nh,kh)).\]

The equation, (\ref{eq-LBE-A}), can be written equivalently as
\begin{equation}\label{eq-LBE}
\left\{ \begin{array}{l}\displaystyle
u^{k+1}_{n+1}-u^k_n=\frac{imh}{2}(v^{k+1}_{n-1}+v^k_n) +\frac{i\alpha h(u^{k+1}_{n+1}+u^k_n)}{2}|v^k_n|^2 +\frac{ih\beta}{2}(v^{k+1}_{n-1}+v^k_n)G(u^k_n,v^k_n), \\  \, \\ \displaystyle v^{k+1}_{n-1}-v^k_n=\frac{imh}{2}(u^{k+1}_{n+1}+u^k_n) +\frac{i\alpha h(v^{k+1}_{n-1}+v^k_n)}{2}|u^k_n|^2 +\frac{ih\beta}{2}(u^{k+1}_{n+1}+u^k_n)G(u^k_n,v^k_n).
\end{array}
\right.
\end{equation}
 Here and in the sequel, we call (\ref{eq-LBE-A}) or (\ref{eq-LBE}) a QLB scheme briefly.

Numerical experiments are given in \cite{L-dellar,succi-benzi,succi2001}  to show the evidence of the convergence of the numerical solutions of (\ref{eq-LBE}) to the nonlinear Dirac equations.  But to our knowledge, there is no rigorous proof of the convergence results on the scheme (\ref{eq-LBE}) or (\ref{eq-LBE-A}). The motivation of this paper is to prove that the solution $(u^{(h)}, v^{(h)})$ given by the scheme (\ref{eq-LBE-A})  (or (\ref{eq-LBE})) is convergent to the strong solutions of (\ref{eq-dirac}) as $h$ goes to $0$.  The main result is stated as follows.

\begin{theorem}\label{thm-main}
 Let $(u_0,v_0)\in L^2(R^1)$ and $\sup\limits_{h\in (0,1)}||(u_0^{(h)},v_0^{(h)})||_{L^2(R^1)}<\infty$ such that
 \begin{equation}\label{eq-mainresult-assum}
 \lim\limits_{h\to 0+} (||u^{(h)}_0-u_0||_{L^2(R^1)}+||v^{(h)}_0-v_0||_{L^2(R^1)})=0.
 \end{equation}
 Then  the QLB scheme (\ref{eq-LBE-A}) with $(u^{(h)}, v^{(h)})(\cdot, t=0)=(u^{(h)}_0, v^{(h)}_0)$ has a unique global solution $(u^{(h)}, v^{(h)})$ for $h\in (0,1)$. Moreover, there holds that
\begin{equation}\label{eq-mainresult}
\lim\limits_{h\to 0+} (||u^{(h)}-u_*||_{C([0,T];L^2(R^1))}+||v^{(h)}-v_*||_{C([0,T];L^2(R^1))})=0.
\end{equation}
for any $T>0$, where $(u_*,v_*)$ is the unique strong solution to (\ref{eq-dirac}) and (\ref{eq-dirac-initialv}).
\end{theorem}

Here the strong solution to (\ref{eq-dirac}) and (\ref{eq-dirac-initialv}) is defined as follows.
\begin{definition}\label{def-weaksolution}
  A pair of functions $(u,v)\in C([0,\infty); L^2(R^1))$ is called a strong solution to (\ref{eq-dirac}) and (\ref{eq-dirac-initialv}) on $R^1\times [0,\infty)$ if there exits a sequence of smooth solutions $(u_{(n)}, v_{(n)}) $ to (\ref{eq-dirac}) on $R^1\times [0,\infty)$ such that
  \[ \lim\limits_{n\to\infty} \big( ||u_{(n)}(\cdot,0)-u_0||_{L^2(R^1)}+||v_{(n)}(\cdot,0)-v_0||_{L^2(R^1)} \big)=0\]
  and
  \[ \lim\limits_{n\to\infty} \big( ||u_{(n)}-u||_{L^2(R^1\times [0,T])}+||v_{(n)}-v||_{L^2(R^1\times [0,T])} \big)=0 \]
   for any $T>0$.
  \end{definition}

The QLB scheme (\ref{eq-LBE-A}) and its equivalent form (\ref{eq-LBE}) are implicit  and nonlinear equations with cubic terms, which bring the difficulties in getting the stability in $L^2$ norms of the solutions. To overcome these difficulties, we make use of their special nonlinear structure and introduce some nonlinear functionals to deal with the nonlinear terms.  More precisely, we first  deduce the explicit estimates (\ref{eq-interaction-1A}) and (\ref{eq-interaction-2A}) for the evolution law from $(|u^k_n|^2, |v^k_n|^2)$ to $(|u^{k+1}_{n+1}|^2,|v^{k+1}_{n-1}|^2)$ and deduce the explicit estimates (\ref{eq-interaction-3B}) and $(\ref{eq-interaction-4B})$  for the evolution law from $(|U^k_n|^2, |V^k_n|^2)$ to $(|U^{k+1}_{n+1}|^2,|V^{k+1}_{n-1}|^2)$ from the implicit homogeneous scheme (\ref{eq-LBE}) and inhomogeneous scheme (\ref{eq-LBE-B}). Here $(U^{k}_{n},V^{k}_{n})$ denotes the difference between the $(n,k)$ components  $(u^k_n,v^k_n)$ and $(\widetilde{u}^k_n,\widetilde{v}^k_n)$ of two solutions to (\ref{eq-LBE}), see section 3. Noticing that  (\ref{eq-interaction-1A}),(\ref{eq-interaction-2A}), (\ref{eq-interaction-3B}) and $(\ref{eq-interaction-4B})$ have quadruple terms with special structures and are analogous to Glimm's estimates for the interactions of waves in \cite{glimm} (see also \cite{bressan}, \cite{dafermos}), we follow the idea from \cite{zhang-zhao} to  introduce a Bony type functional $Q_1$ and a Glimm type functional $F_1(k;\Delta)$, see Definition \ref{def-functional-1} and  Definition \ref{def-functional-2}. Then we can establish the estimates on $F_1$, which  enables us to prove the uniform continuity in  $L^2$ of the solutions $(u^{(h)}, v^{(h)})$ to QLB scheme (\ref{eq-LBE}). And the uniform continuity of the solutions along the characteristic is also proved based on the estimates on the solutions $(u^{(h)}, v^{(h)})$ on characteristics. Such two results imply the relatively compactness of the set of the  solutions $(u^{(h)}, v^{(h)})$, that is, as the mesh size $h$ goes to zero, any sequence of  solutions $(u^{(h)}, v^{(h)})$ has a convergent subsequence in $L^2$. Finally we  estimate the difference between the smooth solution of (\ref{eq-dirac}) and the solutions $(u^{(h)}, v^{(h)})$ by (\ref{eq-LBE}), then prove that every limit of the convergent subsequence of the solutions $(u^{(h)}, v^{(h)})$ is the strong solution of (\ref{eq-dirac}). We remark that Glimm type functional  was  first used by Glimm \cite{glimm} and later by others to establish global existence of small solution to some nonlinear hyperbolic systems,  and that the Bony functional was used  to study the discrete Boltzmann equations, see for instance \cite{bony,bressan,dafermos,ha-tzavaras} and references therein.  There also have been many works on the stability and convergence of the lattice Boltzmann method for other types of partial differential equations, see \cite{junk-yang,junk-yong} for instance and references therein. For the lattice Boltzmann method and its application, see for instance \cite{succi2001}.

The remaining part is organized as follows. In section 2, we establish some point-estimates on the approximate solutions for the scheme. In section 3, we give some local space-time estimates on the differences between two approximate solutions. In section 4 we prove that any sequence of approximate solutions by (\ref{eq-LBE}) has a convergent subsequence in $L^2$. In section 5, we prove that every limit of the convergent sequence of the approximate solutions coincides with  the strong solution of (\ref{eq-dirac}).

\section{Estimates on the solutions to the QLB scheme}\label{section-estimates-1-solu}

\subsection{Homogeneous difference scheme}

We consider the homogeneous scheme (\ref{eq-LBE}) for $h\in (0,1)$ and assume that there exists a constant $C_0>0$ independent of $h$ such that
\begin{equation}\label{assump-1} \sum\limits_{l=-\infty}^{\infty} (|u^{0}_l|^2+|v_l^{0}|^2)h\le C_0. \end{equation}

\begin{lemma}\label{lemma-interaction-0}
For any $h\in (0,1)$, the scheme (\ref{eq-LBE}) is uniquely solvable at each time step. Moreover, for any integers $n$ and $k$ with $k\ge 0$,
 there holds that
\begin{equation}\label{eq-LBE-conserv-1} |u^{k+1}_{n+1}|^2+|v^{k+1}_{n-1}|^2=|u^k_n|^2+|v^k_n|^2 \end{equation}
 and
 \begin{equation}\label{eq-LBE-1}
 \frac{|u_{n+1}^{k+1}|^2-|u^k_n|^2}{h}=\Re\big\{ im(\overline{u^{k+1}_{n+1}+u^k_n})(v^{k+1}_{n-1}+v^k_n)\big\}+e^{k,1}_n,
 \end{equation}
 \begin{equation}\label{eq-LBE-2}
 \frac{|v_{n-1}^{k+1}|^2-|v^k_n|^2}{h}=\Re\big\{ im(\overline{v^{k+1}_{n-1}+v^k_n})(u^{k+1}_{n+1}+u^k_n)\big\}+e^{k,2}_n,
 \end{equation}
 where the remainders are
 \[ e^{k,1}_n= \Re\big\{ i\beta(\overline{u^{k+1}_{n+1}+u^k_n})(v^{k+1}_{n-1}+v^k_n) (\overline{u^k_n}v^k_n+u^k_n\overline{v^k_n})\big\}\]
 and
 \[ e^{k,2}_n=\Re\big\{ i\beta(\overline{v^{k+1}_{n-1}+v^k_n})(u^{k+1}_{n+1}+u^k_n) (\overline{u^k_n}v^k_n+u^k_n\overline{v^k_n})\big\}.\]
 Here and in sequel $\Re z=\frac{z+\overline{z}}{2}$ stands for the real part of $z$ for $z\in C$.
\end{lemma}
{\it Proof.} At time step $t=(k+1)h$, the system (\ref{eq-LBE}) is a linear system for $(u^{k+1}_{n+1},v^{k+1}_{n-1})$  for each pair $(n,k)$. To get the term $(u^{k+1}_{n+1},v^{k+1}_{n-1})$ from the equations (\ref{eq-LBE}), we  compute the determinant $J^k_n$ of coefficients of the term $(u^{k+1}_{n+1},v^{k+1}_{n-1})$ as follows,
\begin{eqnarray*}
J^k_n & =& \det\left(\begin{array}{ll} \qquad 1-\frac{i\alpha}{2}|v^k_n|^2h & -\frac{ih}{2}[m +\beta(\overline{u^k_n}v^k_n+u^k_n\overline{v^k_n})] \\
-\frac{ih}{2}[m +\beta(\overline{u^k_n}v^k_n+u^k_n\overline{v^k_n})] & \qquad 1-\frac{i\alpha}{2}|u^k_n|^2h \end{array}\right) \\
&=& 1-\frac{\alpha^2 h^2 |u^k_n|^2|v^k_n|^2 }{4}+\frac{h^2}{4}[m +\beta(\overline{u^k_n}v^k_n+u^k_n\overline{v^k_n})]^2- i \alpha\frac{|u^k_n|^2+|v^k_n|^2}{2}h^2.
\end{eqnarray*}
Direct computation shows that
\[ |J^k_n|^2 =1+\big(-\frac{\alpha^2h^2|u^k_n|^2|v^k_n|^2}{4}+\frac{(m+G^k_n)^2h^2}{4} \big)^2 +\frac{h^2(m+\beta G^k_n)^2}{2} +\frac{\alpha^2h^2(|u^k_n|^4+|v^k_n|^4)}{4}\ge 1, \]
where $G^k_n=(\overline{u^k_n}v^k_n+u^k_n\overline{v^k_n})\in R^1$.

Therefore, by Cramer's rule, we have unique solution $(u^{k+1}_{n+1},v^{k+1}_{n-1})$ for the equations (\ref{eq-LBE})  and prove the solvability of the equations (\ref{eq-LBE}).

Now multiplying the first and second equations in (\ref{eq-LBE}) by $\overline{u^{k+1}_{n+1}+u^k_n}$ and $\overline{v^{k+1}_{n-1}+v^k_n}$ respectively and taking their real parts, we can have (\ref{eq-LBE-1}) and (\ref{eq-LBE-2}).

Finally, taking the sum of (\ref{eq-LBE-1}) and (\ref{eq-LBE-2}) gives (\ref{eq-LBE-conserv-1}).
The proof is complete. $\Box$.

Due to Lemma \ref{lemma-interaction-0}, the scheme (\ref{eq-LBE}) has a global solution. Let $\{(u^k_n, v^k_n)\}$ be the solution to (\ref{eq-LBE}) in the sequel, and we have the following.
\begin{corollary}\label{coro-conservation}
For any integer $k\ge 0$, there holds that
\[ \sum_{n=-\infty}^{\infty}(|u^k_n|^2+|v^k_n|^2)=\sum_{n=-\infty}^{\infty}(|u^0_n|^2+|v^0_n|^2).\]
\end{corollary}
{\it Proof.} Taking the sum of (\ref{eq-LBE-conserv-1}) over $n$ yields that
\[ \sum_{n=-\infty}^{\infty}(|u^{k+1}_n|^2+|v^{k+1}_n|^2)=\sum_{n=-\infty}^{\infty}(|u^k_n|^2+|v^k_n|^2), \]
which gives the desired result by induction on $n$ and completes the  proof.$\Box$

We consider the scheme (\ref{eq-LBE}) on the triangle domains. For any integers $n_1, k_1$ and $k_0$ with $0\le k_0\le k_1$, denote
\[ \Delta(n_1,k_1; k_0)=\{ (n,k)\big| \mbox{$n,k$ are integers and}\, n_1-k_1+k\le n\le n_1+k_1-k,\, k_0\le k\le k_1\},\] see Fig. \ref{fig-domain1}.
\begin{figure}[h]
\begin{center}
\unitlength=10mm
\begin{picture}(10,3.5)
\thicklines
\put(0,0){\line(1,0){10}}
\put(1,0){\line(5,4){4}}
\put(9,0){\line(-5,4){4}}
\put(0,-0.5){$(n_1-k_1+k_0,k_0)$}
\put(8,-0.5){$(n_1+k_1-k_0,k_0)$}\put(10.3,0){$k=k_0$}
\put(4.5,3.5){$(n_1,k_1)$}
\put(4,1.5){$\Delta(n_1,k_1; k_0)$}
\end{picture}
\caption{The set $\Delta(n_1,k_1; k_0)$}\label{fig-domain1}
\end{center}
\end{figure}

Taking the summation of (\ref{eq-LBE-conserv-1}) over $\Delta(n_1, k_1; k)$ gives the following.
\begin{lemma}\label{lemma-sum-1}
For $0\le k_0+1\le k\le k_1$ and $-\infty<n_1<\infty$, there holds
\[ \sum\limits_{k_1-k\le j\le k_1-k_0} |v_{n_1-1-j}^{k_1+1-j}|^2+\sum\limits_{k_1-k\le j\le k_1-k_0} |u_{n_1+1+j}^{k_1+1-j}|^2 \le\sum\limits_{l=n_1-k_1+k_0}^{n_1+k_1-k_0} (|u^{k_0}_l|^2+|v_l^{k_0}|^2)\]
and
\[ \sum\limits_{l=n_1-k_1+k+1}^{n_1+k_1-k+1} |u^{k+1}_l|^2 +\sum\limits_{l=n_1-k_1+k-1}^{n_1+k_1-k-1} |v^{k+1}_l|^2  \le \sum\limits_{l=n_1-k_1+k_0}^{n_1+k_1-k_0} (|u^{k_0}_l|^2+|v_l^{k_0}|^2). \]
Therefore,
\[ \sum\limits_{0\le j\le k_1} |v_{n_1-1-j}^{k_1+1-j}|^2+\sum\limits_{0\le j\le k_1} |u_{n_1+1+j}^{k_1+1-j}|^2 \le\sum\limits_{l=-\infty}^{\infty} (|u^{k_0}_l|^2+|v_l^{k_0}|^2).\]
\end{lemma}
{\it Proof.}
By (\ref{eq-LBE-conserv-1}), we have
\[ \sum_{(n,j)\in \Delta(n_1,k_1; k)}( |u^{j+1}_{n+1}|^2+|v^{j+1}_{n-1}|^2-|u^j_n|^2-|v^j_n|^2)=0,\]
where the cancelation of terms for $(n,j)$ in the interior of $\Delta(n_1,k_1; k)$  gives the proof of the lemma. The proof is complete.$\Box$

Then, we have the pointwise estimates as follows.
\begin{lemma}\label{lemma-pointwise-1}
There exist a constant $C_1>0$, independent of $h$ and $(n,k)$, such that
\begin{equation}
|u^{k+1}_{n+1}|\le C_1|u^0_{n-k}|+C_1\sqrt{kh} \label{eq-pointwise-1}
\end{equation}
and
\begin{equation}
|v^{k+1}_{n-1}|\le C_1|v^0_{n+k}|+C_1\sqrt{kh} \label{eq-pointwise-2}
\end{equation}
for $k\ge 0$ and $-\infty<n<\infty$.
\end{lemma}
{\it Proof.} For $0\le j\le k$,
(\ref{eq-LBE}) gives that
\begin{eqnarray*}
|(1-\frac{i\alpha h|v^{k-j}_{n-j}|^2}{2})||u^{k+1-j}_{n+1-j}| &\le & |u^{k-j}_{n-j}|+\frac{mh}{2}(|v^{k+1-j}_{n-1-j}|+|v^{k-j}_{n-j}|)
\\ &\, & +h|\beta||u^{k-j}_{n-j}||v^{k-j}_{n-j}| (|v^{k+1-j}_{n-1-j}|+|v^{k-j}_{n-j}|)\\
&\le & |u^{k-j}_{n-j}|\exp\{4|\beta|h(|v^{k+1-j}_{n-1-j}|^2+|v^{k-j}_{n-j}|^2)\}  \\ &\,& +\frac{mh}{2}(|v^{k+1-j}_{n-1-j}|+|v^{k-j}_{n-j}|).
\end{eqnarray*}
Then
\[|u^{k+1-j}_{n+1-j}|\le |u^{k-j}_{n-j}|\exp\{4|\beta|h(|v^{k+1-j}_{n-1-j}|^2+|v^{k-j}_{n-j}|^2)\} +\frac{mh}{2}(|v^{k+1-j}_{n-1-j}|+|v^{k-j}_{n-j}|),\]
which leads to the following,
\begin{eqnarray*}
|u^{k+1}_{n+1}| &\le & \{|u^0_{n-k}|+mh \sum\limits_{0\le j\le k}(|v^{k+1-j}_{n-1-j}|+|v^{k-j}_{n-j}|)\}\exp\{4|\beta|h\sum\limits_{0\le j\le k}(|v^{k+1-j}_{n-1-j}|^2+|v^{k-j}_{n-j}|^2)\} \\
&\le & \big\{|u^0_{n-k}|+mh \sqrt{\sum\limits_{0\le j\le k}4(|v^{k+1-j}_{n-1-j}|^2+|v^{k-j}_{n-j}|^2)} \sqrt{\sum\limits_{0\le j\le k} 1}\big\}\exp(8|\beta|C_0)\\
&\le & \{|u^0_{n-k}|+4m\sqrt{C_0}\sqrt{kh}\}\exp(8|\beta|C_0),
\end{eqnarray*}
where we use Lemma \ref{lemma-sum-1} and the assumption (\ref{assump-1}) to get last two inequalities. Therefore (\ref{eq-pointwise-1}) is proved.

The inequality (\ref{eq-pointwise-2})  for $v^{k+1}_{n-1}$ can be proved in the same way. Thus, the proof is complete. $\Box$

As one of its consequence, there holds the following.
\begin{lemma}\label{lemma-product-0}
Let $T\in [0,\infty)$. If $0<k_0 \le k_1\le T/h$, then
\begin{eqnarray*}
  \sum^{k_1}_{k=k_0} \sum^{\infty}_{-\infty} |u^k_n|^2|v^k_n|^2 h^2
  \le 4C_1^2 \sum^{\infty}_{n=-\infty} \big(|u^0_n|^2h \sum_{l=n+k_0}^{n+k_1}|v^0_l|^2h\big) \\
  +(4C_1^4+2C_1^2)T(k_1-k_0)h \sum^{\infty}_{-\infty} (|u^0_n|^2+|v^0_n|^2)h.
\end{eqnarray*}
Therefore
\begin{eqnarray*}
\displaystyle  \sum^{k_1}_{k=0} \sum^{\infty}_{-\infty} |u^k_n|^2|v^k_n|^2 h^2
\displaystyle  \le  4C_1^2C_0^2 +(4C_1^4+2C_1^2)C_0T^2.
\end{eqnarray*}
Here  $C_0$ and $C_1$ are the constants given by (\ref{assump-1}) and by Lemma \ref{lemma-pointwise-1}.
\end{lemma}
{\it Proof.} Let
$\displaystyle D_0(k_0, k_1)=\sum^{k_1}_{k=k_0} \sum^{\infty}_{-\infty} |u^k_n|^2|v^k_n|^2$. Then by Lemma \ref{lemma-pointwise-1}, we have
\begin{eqnarray*}
D_0(k_0, k_1)&\le& \sum^{k_1}_{k=k_0} \sum^{\infty}_{-\infty}2C_1( |u^0_{n-k}|^2+kh)|v^k_n|^2 \\
&\le& \sum^{k_1}_{k=k_0} \sum^{\infty}_{-\infty}2C_1|u^0_{n-k}|^2|v^k_n|^2 +2C_1T\sum^{k_1}_{k=k_0}
\sum^{\infty}_{-\infty}|v^k_n|^2 \\
&\le& \sum^{k_1}_{k=k_0} \sum^{\infty}_{-\infty}4C_1^2|u^0_{n-k}|^2(|v^0_{n+k}|^2+kh) +2C_1T\sum^{k_1}_{k=k_0}
\sum^{\infty}_{-\infty}|v^k_n|^2 \\
&\le& \sum^{k_1}_{k=k_0} \sum^{\infty}_{-\infty}4C_1^2|u^0_{n-k}|^2|v^0_{n+k}|^2 +(4C_1^2+2C_1)T\sum^{k_1}_{k=k_0}
\sum^{\infty}_{-\infty}(|u^k_n|^2+|v^k_n|^2),
\end{eqnarray*}
where we use the following,
\begin{eqnarray*}
\sum^{k_1}_{k=k_0} \sum^{\infty}_{-\infty}|u^0_{n-k}|^2|v^0_{n+k}|^2=\sum^{\infty}_{n=-\infty} \big(|u^0_n|^2 \sum_{l=n+k_0}^{n+k_1}|v^0_l|^2\big)
\end{eqnarray*}
and
\begin{eqnarray*}
\sum^{\infty}_{-\infty}(|u^k_n|^2+|v^k_n|^2)=\sum^{\infty}_{-\infty}(|u^0_n|^2+|v^0_n|^2).
\end{eqnarray*}
Therefore we can conclude the result and  the proof is complete.$\Box$

Now we consider the evolution laws for $(|u^k_n|^2, |v^k_n|^2)$. At first, we  deal with the remainders  $e^{k,1}_n$ and $e^{k,2}_n$ given by Lemma \ref{lemma-interaction-0}.

Direct computation by (\ref{eq-LBE-1}) and (\ref{eq-LBE-2}) shows the following.
\begin{lemma}\label{lemma-interaction-1}
There holds that
\[ |e^{k,1}_n|\le |\beta|e^k_n, \quad |e^{k,2}_n|\le |\beta|e^k_n, \]
for $k\ge 0$ and $-\infty<n<\infty$, where
\[ e^k_n=(|u^{k+1}_{n+1}|^2+|u^k_n|^2)|v^k_n|^2+(|v^{k+1}_{n+1}|^2+|v^k_n|^2)|u^k_n|^2.\]
Therefore,
\begin{equation}\label{eq-interaction-1}
\Big|\frac{|u_{n+1}^{k+1}|^2-|u^k_n|^2}{h}\Big|\le 8m(|u^k_n|^2+|v^k_n|^2)+|\beta|e^k_n
\end{equation}
and
\begin{equation}\label{eq-interaction-2}
\Big|\frac{|v_{n-1}^{k+1}|^2-|v^k_n|^2}{h}\Big|\le 8m(|u^k_n|^2+|v^k_n|^2)+|\beta|e^k_n.
\end{equation}
\end{lemma}

This lemma enable us to get the local  estimates on the evolution of $(|u^{(h)}|^2, |v^{(h)}|^2)$ as follows.
\begin{lemma}\label{lemma-interaction-2}
If $(|u^k_n|^2+|v^k_n|^2)h\le \min\{\frac{1}{4|\beta|},\frac{1}{2}\}$, then
\[ 0\le e^k_n\le 8|u^k_n|^2|v^k_n|^2+8m(|u^k_n|^2+|v^k_n|^2).\]
Therefore,
\begin{equation}\label{eq-interaction-1A}
\Big|\frac{|u_{n+1}^{k+1}|^2-|u^k_n|^2}{h}\Big|\le C_{\beta}(|u^k_n|^2+|v^k_n|^2)+C_{\beta}|u^k_n|^2|v^k_n|^2
\end{equation}
and
\begin{equation}\label{eq-interaction-2A}
\Big|\frac{|v_{n-1}^{k+1}|^2-|v^k_n|^2}{h}\Big|\le C_{\beta}(|u^k_n|^2+|v^k_n|^2)+C_{\beta}|u^k_n|^2|v^k_n|^2.
\end{equation}
Here $C_{\beta}=8m+16|\beta|m +16|\beta|$.
\end{lemma}
{\it Proof.} By Lemma \ref{lemma-interaction-1},
\[ |u^{k+1}_{n+1}|^2\le |u^k_n|^2+8mh(|u^k_n|^2+|v^k_n|^2)+2|\beta| h e^k_n\] and
\[ |v^{k+1}_{n-1}|^2\le |v^k_n|^2+8mh(|u^k_n|^2+|v^k_n|^2)+2|\beta| h e^k_n,\]
which leads to the following,
\begin{eqnarray}
e^k_n &=& (|u^k_n|^2+|u^{k+1}_{n+1}|^2)|v^k_n|^2+(|v^k_n|^2+|v^{k+1}_{n-1}|^2)|u^k_n|^2 \nonumber\\
 &\le& 4|u^k_n|^2|v^k_n|^2+ 8mh(|u^k_n|^2+|v^k_n|^2)^2 +2|\beta|e^k_n (|u^k_n|^2+|v^k_n|^2)h. \label{eq-error1-1}
\end{eqnarray}
For $\beta \neq 0$,  we can get the estimate for $e^k_n$ from (\ref{eq-error1-1}) for $(|u^k_n|^2+|v^k_n|^2)h\le \min\{\frac{1}{4|\beta|},\frac{1}{2}\}$; while for $\beta = 0$, we can get the estimate for $e^k_n$ from (\ref{eq-error1-1}) for $(|u^k_n|^2+|v^k_n|^2)h\le \frac{1}{2}$.

Moreover, plugging the estimate on $e^k_n$ into (\ref{eq-interaction-1}) and (\ref{eq-interaction-2}) gives (\ref{eq-interaction-1A}) and (\ref{eq-interaction-2A}). The proof is complete.$\Box$

\subsection{Inhomogeneous difference scheme}

Let $(\widetilde{u}^{(h)},\widetilde{v}^{(h)})$ be the solution to the following scheme
\begin{equation}\label{eq-LBE-B}
\left\{ \begin{array}{l}
 \displaystyle\frac{\widetilde{u}^{k+1}_{n+1}-\widetilde{u}^k_n}{h}
 \displaystyle =\frac{im}{2}(\widetilde{v}^{k+1}_{n-1}+\widetilde{v}^k_n)
 \displaystyle +\frac{i\alpha (\widetilde{u}^{k+1}_{n+1}+\widetilde{u}^k_n)}{2}|v^k_n|^2
 \displaystyle+\frac{i\beta}{2}(\widetilde{v}^{k+1}_{n-1}+\widetilde{v}^k_n)\widetilde{G}^k_n +g^{k,1}_n,
 \\  \, \\  \displaystyle \frac{\widetilde{v}^{k+1}_{n-1}-\widetilde{v}^k_n}{h} \displaystyle=\frac{im}{2}(\widetilde{u}^{k+1}_{n+1}+\widetilde{u}^k_n)
 \displaystyle +\frac{i\alpha (\widetilde{v}^{k+1}_{n-1}+\widetilde{v}^k_n)}{2}|\widetilde{u}^k_n|^2 \displaystyle+\frac{i\beta}{2}(\widetilde{u}^{k+1}_{n+1}+\widetilde{u}^k_n)\widetilde{G}^k_n +g^{k,2}_n,
\end{array}
\right.
\end{equation}
for given data $g^{k,1}_n$ and $g^{k,2}_n$ with integers $k$ and $n$ satisfying $k\ge 0$ and $-\infty<n<\infty$. Here
\[\widetilde{G}^k_n=G(\widetilde{u}^k_n,\widetilde{v}^k_n)=\overline{\widetilde{u}^k_n}\widetilde{v}^k_n + \widetilde{u}^k_n\overline{\widetilde{v}^k_n}\]
and
the function $(\widetilde{u}^{(h)},\widetilde{v}^{(h)}) $ is piecewise-constants valued, and  satisfies
\begin{equation}\label{eq-funct-value}
(\widetilde{u}^{(h)}(x,t),\widetilde{v}^{(h)}(x,t))=(\widetilde{u}^k_n,\widetilde{v}^k_n), \, \, (x,t)\in [nh,(n+1)h)\times [kh,(k+1)h)
\end{equation}
for the any integers $n$ and $k\ge 0$, where
\[(\widetilde{u}^k_n,\widetilde{v}^k_n)=(\widetilde{u}^{(h)}(nh,kh),\widetilde{v}^{(h)}(nh,kh)).\]

As in the proof of Lemma \ref{lemma-interaction-1} for homogeneous case (\ref{eq-LBE}), we carry out the same argument to derive the following.
\begin{lemma}\label{lemma-interaction-1-B}
For any $h\in (0,1)$, the scheme (\ref{eq-LBE-B}) is uniquely solvable at each time step. Moreover, for any integers $n$ and $k$ with $k\ge 0$,
 there holds that
\begin{equation}\label{eq-interaction-1-B}
 \Big|\frac{|\widetilde{u}_{n+1}^{k+1}|^2-|\widetilde{u}^k_n|^2}{h}\Big|\le
  (8m+2)(|\widetilde{u}_{n+1}^{k+1}|^2+|\widetilde{v}_{n-1}^{k+1}|^2
  +|\widetilde{u}^k_n|^2+|\widetilde{v}^k_n|^2)+|\beta|\widetilde{e}^k_n+ |g^k_n|^2
\end{equation}
and
\begin{equation}\label{eq-interaction-2-B}
\Big|\frac{|\widetilde{v}_{n-1}^{k+1}|^2-|\widetilde{v}^k_n|^2}{h}\Big|\le (8m+2)(|\widetilde{u}_{n+1}^{k+1}|^2+|\widetilde{v}_{n-1}^{k+1}|^2 +|\widetilde{u}^k_n|^2+|\widetilde{v}^k_n|^2) +|\beta|e^k_n + |g^k_n|^2,
\end{equation}
where
 \[ \widetilde{e}^k_n=(|\widetilde{u}^{k+1}_{n+1}|^2+|\widetilde{u}^k_n|^2)|\widetilde{v}^k_n|^2 +(|\widetilde{v}^{k+1}_{n-1}|^2+|\widetilde{v}^k_n|^2)|\widetilde{u}^k_n|^2\]
 and
 \[ |g^k_n|^2=|g^{k,1}_n|^2+|g^{k,2}_n|^2.\]
\end{lemma}

Then we have the following evolution estimates for $(\widetilde{u}^{(h)},\widetilde{v}^{(h)})$.
\begin{lemma}\label{lemma-interaction-2-B}
There exist constants $\delta_1>0$ and $C_2>0$ such that if $h\in (0,\frac{1}{2})$ and if $(|\widetilde{u}^k_n|^2+|\widetilde{v}^k_n|^2)h\le \delta_1$ then
\begin{equation}\label{eq-interaction-1A-B}
\Big|\frac{|\widetilde{u}_{n+1}^{k+1}|^2-|\widetilde{u}^k_n|^2}{h}\Big|\le C_{2}\Big((|\widetilde{u}^k_n|^2+|\widetilde{v}^k_n|^2)+|\widetilde{u}^k_n|^2|\widetilde{v}^k_n|^2 +|g^k_n|^2 \Big)
\end{equation}
and
\begin{equation}\label{eq-interaction-2A-B}
\Big|\frac{|\widetilde{v}_{n-1}^{k+1}|^2-|\widetilde{v}^k_n|^2}{h}\Big|\le C_{2}\Big((|\widetilde{u}^k_n|^2+|\widetilde{v}^k_n|^2)+|\widetilde{u}^k_n|^2|\widetilde{v}^k_n|^2 +|g^k_n|^2 \Big).
\end{equation}
\end{lemma}
{\it Proof.} At first, as in the proof of Lemma \ref{lemma-interaction-0}, we multiply the first equation in (\ref{eq-LBE-B}) by $ \overline{\widetilde{u}^{k+1}_{n+1}+\widetilde{u}^k_n}$ and the second equation in
(\ref{eq-LBE-B}) by $ \overline{\widetilde{v}^{k+1}_{n-1}+\widetilde{v}^k_n}$, and take the sum of their real parts to deduce that
\[ |\widetilde{u}^{k+1}_{n+1}|^2+|\widetilde{v}^{k+1}_{n-1}|^2=|\widetilde{u}^k_n|^2+|\widetilde{v}^k_n|^2+h \Re \{ g^{k,1}_n(\overline{\widetilde{u}^{k+1}_{n+1}+\widetilde{u}^k_n})+ g^{k,2}_n(\overline{\widetilde{v}^{k+1}_{n-1}+\widetilde{v}^k_n})\}.\]
Then,
\[ |\widetilde{u}^{k+1}_{n+1}|^2+|\widetilde{v}^{k+1}_{n-1}|^2
\le |\widetilde{u}^k_n|^2+|\widetilde{v}^k_n|^2+h\{ |g^k_n|^2+|\widetilde{u}^{k+1}_{n+1}|^2+|\widetilde{v}^{k+1}_{n-1}|^2 +|\widetilde{u}^k_n|^2+|\widetilde{v}^k_n|^2\},
\]
which gives  for $h\in (0, \frac{1}{2})$ that
\begin{equation}\label{eq-interaction-2A-Bpf}
|\widetilde{u}^{k+1}_{n+1}|^2+|\widetilde{v}^{k+1}_{n-1}|^2
\le 2(|\widetilde{u}^k_n|^2+|\widetilde{v}^k_n|^2)+ 2|g^k_n|^2.
\end{equation}
Plugging (\ref{eq-interaction-2A-Bpf}) into (\ref{eq-interaction-1-B}) and (\ref{eq-interaction-2-B}) yields that
\begin{equation}\label{eq-interaction-1-Bpf}
\Big|\frac{|\widetilde{u}_{n+1}^{k+1}|^2-|\widetilde{u}^k_n|^2}{h}\Big|\le 4(8m+2)( |\widetilde{u}^k_n|^2+|\widetilde{v}^k_n|^2)+|\beta|\widetilde{e}^k_n+(16m+5) |g^k_n|^2
\end{equation}
and
\begin{equation}\label{eq-interaction-2-Bpf}
\Big|\frac{|\widetilde{v}_{n-1}^{k+1}|^2-|\widetilde{v}^k_n|^2}{h}\Big|\le 4(8m+2)( |\widetilde{u}^k_n|^2+|\widetilde{v}^k_n|^2) +|\beta|e^k_n +(16m+5) |g^k_n|^2,
\end{equation}
which enable us to carry out same argument as in the proof of Lemma \ref{lemma-interaction-2} to give  (\ref{eq-interaction-1A-B}) and (\ref{eq-interaction-2A-B}). The proof is complete.$\Box$

\section{$L^2$- stability Estimates on the QLB schemes}

\subsection{Estimates on the difference of solutions}

Let $\{(u^k_n,v^k_n)\}_{k,n}$ be given by scheme (\ref{eq-LBE}) and $\{(\widetilde{u}^k_n,\widetilde{v}^k_n)\}_{k,n}$ be given by (\ref{eq-LBE-B}).  Denote
\[ U^k_n=\widetilde{u}^k_n-u^k_n, \quad V^k_n=\widetilde{v}^k_n-v^k_n\]
for integers $k$ and $n$ with $k\ge 0$ and $-\infty<n<\infty$.

Then
\begin{equation}\label{eq-LBE-3}
\frac{U^{k+1}_{n+1}-U^k_n}{h}=im \frac{V^{k+1}_{n-1}+V^k_n}{2} +i\alpha q^{k,1}_n+i\beta q^{k,2}_n+ g^{k,1}_n
\end{equation}
and
\begin{equation}\label{eq-LBE-4}
\frac{V^{k+1}_{n-1}-V^k_n}{h}=im \frac{U^{k+1}_{n+1}+U^k_n}{2} +i\alpha q^{k,3}_n+i\beta q^{k,4}_n +g^{k,2}_n,
\end{equation}
where
\[ q^{k,1}_n=\frac{U^{k+1}_{n+1}+U^k_n}{2}|\widetilde{v}^k_n|^2 +\frac{u^{k+1}_{n+1}+u^k_n}{2}(V^k_n\overline{\widetilde{v}^k_n}+v^k_n\overline{V^k_n}),\]
\[ q^{k,2}_n=\frac{1}{2}\Big\{(V^{k+1}_{n-1}+V^k_n)G(\widetilde{u}^k_n, \widetilde{v}^k_n) +(v^{k+1}_{n-1}+v^k_n) G(U^k_n, \widetilde{v}^k_n) +(v^{k+1}_{n-1}+v^k_n) G(u^k_n, V^k_n)\Big\},
\]
\[ q^{k,3}_n=\frac{V^{k+1}_{n-1}+V^k_n}{2}|\widetilde{u}^k_n|^2 +\frac{v^{k+1}_{n-1}+v^k_n}{2}(U^k_n\overline{\widetilde{u}^k_n}+u^k_n\overline{U^k_n})\]
and
\[ q^{k,4}_n=\frac{1}{2}\Big\{(U^{k+1}_{n+1}+U^k_n)G(\widetilde{u}^k_n, \widetilde{v}^k_n) +(u^{k+1}_{n+1}+u^k_n) G(U^k_n, \widetilde{v}^k_n) +(u^{k+1}_{n+1}+u^k_n) G(u^k_n, V^k_n)\Big\}.
\]

Direct computation leads to the following estimates on $(U^{k+1}_{n+1}, V^{k+1}_{n-1})$.
\begin{lemma}\label{lemma-interaction-3}
There exists a constant $\widehat{C}_3>0$ such that
\begin{equation}\label{eq-interaction-3A}
|\frac{|U^{k+1}_{n+1}|^2-|U^k_n|^2}{h}|\le \widehat{C}_3(|U^{k+1}_{n+1}|^2+|V^{k+1}_{n-1}|^2+|U^k_n|^2+|V^k_n|^2+\widehat{E}^k_n+|g^k_n|^2)
\end{equation}
and
\begin{equation}\label{eq-interaction-4A}
|\frac{|V^{k+1}_{n-1}|^2-|V^k_n|^2}{h}|\le \widehat{C}_3(|U^{k+1}_{n+1}|^2+|V^{k+1}_{n-1}|^2+|U^k_n|^2+|V^k_n|^2+\widehat{E}^k_n+|g^k_n|^2)
\end{equation}
for any $k\ge 0$, $-\infty<n<\infty$ and $h\in (0,\frac{1}{4})$,
 where
\begin{eqnarray*} \widehat{E}^k_n &=&(|U^{k+1}_{n+1}|^2+|U^k_n|^2) (|v^{k+1}_{n-1}|^2+|v^k_n|^2+|\widetilde{v}^{k+1}_{n-1}|^2+|\widetilde{v}^k_n|^2) \\ &+&(|V^{k+1}_{n-1}|^2+|V^k_n|^2) (|u^{k+1}_{n+1}|^2+|u^k_n|^2+|\widetilde{u}^{k+1}_{n+1}|^2+|\widetilde{u}^k_n|^2).
\end{eqnarray*}
\end{lemma}
{\it Proof.} Indeed, multiplying the equations (\ref{eq-LBE-3}) and (\ref{eq-LBE-4}) by $\overline{U^{k+1}_{n+1}+U^k_n}$ and $\overline{V^{k+1}_{n-1}+V^k_n}$ respectively and taking the real parts, we have
\begin{equation}\label{eq-interaction-3}
\frac{|U^{k+1}_{n+1}|^2-|U^k_n|^2}{h}=\Re \Big\{ im(\overline{U^{k+1}_{n+1}+U^k_n})(V^{k+1}_{n-1}+V^k_n)\Big\}+E^{k,1}_n
\end{equation}
and
\begin{equation}\label{eq-interaction-4}
\frac{|V^{k+1}_{n-1}|^2-|V^k_n|^2}{h}=\Re \Big\{ im(\overline{V^{k+1}_{n-1}+V^k_n})(U^{k+1}_{n+1}+U^k_n)\Big\}+E^{k,2}_n,
\end{equation}
where
\[ E^{k,1}_n=\Re\Big\{ (\overline{U^{k+1}_{n+1}+U^k_n})(i\alpha q^{k,1}_n +i\beta q^{k,2}_n +g^{k,1}_n ) \Big\}\]
and
\[ E^{k,2}_n=\Re\Big\{ (\overline{V^{k+1}_{n-1}+V^k_n})(i\alpha q^{k,3}_n+ i\beta q^{k,4}_n  +g^{k,2}_n\Big\}.\]

 Then applying Cauchy -Schwarz inequality to the righthand sides in (\ref{eq-interaction-3}) and (\ref{eq-interaction-4}) respectively leads to (\ref{eq-interaction-3A}) and
 (\ref{eq-interaction-4A}). Thus the proof is complete. $\Box$

In addition, the above estimates could be modified to the more exact ones as follow.
\begin{lemma}\label{lemma-interaction-4}
There exist constants $\delta_2>0$ and $C_3>0$ such that if $h\in (0, h_1)$ and if $(|u^k_n|^2+|v^k_n|^2)h\le \delta_2$ and  $(|\widetilde{u}^k_n|^2+|\widetilde{v}^k_n|^2)h\le \delta_2$ and if $|g^k_n|^2 h\le \delta_2$, then
\begin{equation}\label{eq-interaction-3B}
|\frac{|U^{k+1}_{n+1}|^2-|U^k_n|^2}{h}|\le C_3(|U^k_n|^2+|V^k_n|^2+ E^k_n+|g^k_n|^2)
\end{equation}
and
\begin{equation}\label{eq-interaction-4B}
|\frac{|V^{k+1}_{n-1}|^2-|V^k_n|^2}{h}|\le C_3(|U^k_n|^2+|V^k_n|^2+ E^k_n+|g^k_n|^2),
\end{equation}
where
\[ E^k_n=|U^k_n|^2(|v^k_n|^2+|\widetilde{v}^k_n|^2)+|V^k_n|^2(|u^k_n|^2+|\widetilde{u}^k_n|^2).\]
Here $h_1=\min\{\frac{1}{2\widehat{C}_3},\frac{1}{4}\}$.
\end{lemma}
{\it Proof.}
At first  Lemma \ref{lemma-interaction-3} gives the following,
\[ |U^{k+1}_{n+1}|^2+|V^{k+1}_{n-1}|^2 \le \widehat{C}_3(h+1)(|U^k_n|^2+|V^k_n|^2)+\widehat{C}_3h (|U^{k+1}_{n+1}|^2 +|V^{k+1}_{n-1}|^2 + \widehat{E}^k_n+|g^k_n|^2),\]
which implies that for $h\in (0, \frac{1}{2\widehat{C}_3})$ there holds that
\[ |U^{k+1}_{n+1}|^2+|V^{k+1}_{n-1}|^2 \le 2\widehat{C}_3(h+1)(|U^k_n|^2+|V^k_n|^2)+2\widehat{C}_3h ( \widehat{E}^k_n+|g^k_n|^2).\]
Therefore, plugging the above inequality into (\ref{eq-interaction-3A}) and
 (\ref{eq-interaction-4A}),  we get
\begin{equation}\label{eq-interaction-3B-0}
|\frac{|U^{k+1}_{n+1}|^2-|U^k_n|^2}{h}|\le C_3^{\prime}(|U^k_n|^2+|V^k_n|^2+ \widehat{E}^k_n+|g^k_n|^2)
\end{equation}
and
\begin{equation}\label{eq-interaction-4B-0}
|\frac{|V^{k+1}_{n-1}|^2-|V^k_n|^2}{h}|\le C_3^{\prime}(|U^k_n|^2+|V^k_n|^2+ \widehat{E}^k_n+|g^k_n|^2)
\end{equation}
for $h\in (0, h_1)$ with some constant $C^{\prime}_3>0$ depending only on $\widehat{C}_3$.

Next, we assume that $h\in (0,h_1)$ and assume that $(|u^k_n|^2+|v^k_n|^2)h\le \delta_2$ and   $(|\widetilde{u}^k_n|^2+|\widetilde{v}^k_n|^2)h\le \delta_2$ and $|g^k_n|^2h\le \delta_2$, where $\delta_2\in (0,\frac{1}{(C_1+C_{\beta})+1})$ is a constant to be specified later.

 By (\ref{eq-interaction-3B-0}) and (\ref{eq-interaction-4B-0}), we have
\begin{equation}\label{eq-interaction1-pf}
|U^{k+1}_{n+1}|^2+|U^k_n|^2\le 2|U^k_n|^2+C_3^{\prime} h(|U^k_n|^2+|V^k_n|^2+\widehat{E}^k_n+|g^k_n|^2)
\end{equation}
and
\begin{equation}\label{eq-interaction2-pf}
|V^{k+1}_{n-1}|^2+|V^k_n|^2\le 2|V^k_n|^2+C_3^{\prime}h(|U^k_n|^2+|V^k_n|^2+\widehat{E}^k_n+|g^k_n|^2);
\end{equation}
and by Lemma \ref{lemma-interaction-2} and Lemma \ref{lemma-interaction-2-B} we have
\begin{eqnarray}
|u^{k+1}_{n+1}|^2+|u^k_n|^2 +|\widetilde{u}^{k+1}_{n+1}|^2+|\widetilde{u}^k_n|^2  &\le & 2|u^k_n|^2+2|\widetilde{u}^k_n|^2 + C^{\prime}_2 h(|u^k_n|^2|v^k_n|^2+|\widetilde{u}^k_n|^2|\widetilde{v}^k_n|^2 )
 \nonumber\\
&+ & C^{\prime}_2h(|u^k_n|^2+|v^k_n|^2+ |\widetilde{u}^k_n|^2+|\widetilde{v}^k_n|^2+|g^k_n|^2)\nonumber\\
&\le& 3(|u^k_n|^2+|\widetilde{u}^k_n|^2)+3\delta_2 \label{eq-interaction3-pf}
\end{eqnarray}
and
\begin{eqnarray} \label{eq-interaction4-pf}
|v^{k+1}_{n-1}|^2+|v^k_n|^2 +|\widetilde{v}^{k+1}_{n-1}|^2+|\widetilde{v}^k_n|^2
 \le 3(|v^k_n|^2+|\widetilde{v}^k_n|^2)+3\delta_2,
\end{eqnarray}
where $C^{\prime}_2=C_1+C_{\beta}$.

Then, (\ref{eq-interaction1-pf}-\ref{eq-interaction4-pf}) gives the following,
\begin{eqnarray}
\widehat{E}^k_n &\le& \big\{ 2|U^k_n|^2+C^{\prime}_3h(|U^k_n|^2+|V^k_n|^2+\widehat{E}^k_n+|g^k_n|^2)\big\} \big\{ 3(|v^k_n|^2+|\widetilde{v}^k_n|^2)+3\delta_2 \big\} \nonumber\\
&\, & +\big\{2|V^k_n|^2+C^{\prime}_3h(|U^k_n|^2+|V^k_n|^2+\widehat{E}^k_n+|g^k_n|^2)   \big\}\big\{3(|u^k_n|^2+|\widetilde{u}^k_n|^2)+3\delta_2 \big\} \nonumber\\
&\le& 6E^k_n+ C^{\prime\prime}_3\delta_2(|U^k_n|^2+|V^k_n|^2+|g^k_n|^2)+ C^{\prime\prime\prime}_3\delta_2 \widehat{E}^k_n \label{eq-remainder-1}
\end{eqnarray}
for $C^{\prime\prime}_3=12+6C^{\prime}_3(2+h)$  and $C^{\prime\prime\prime}_3=6C^{\prime}_3(2+h)$.

Now we choose a suitable constant $\delta_2>0$ such that $C^{\prime\prime\prime}_3\delta_2<\frac{1}{2}$. Then (\ref{eq-remainder-1}) gives the estimate
 \[\widehat{E}^k_n\le C^{*}_3 (|U^k_n|^2+|V^k_n|^2+ E^k_n+|g^k_n|^2)\] for a positive constant $C^{*}_3=(12+4C^{\prime\prime}_3\delta_2)$, which together with (\ref{eq-interaction-3B-0}) and (\ref{eq-interaction-4B-0})  leads to the desired estimates (\ref{eq-interaction-3B}) and (\ref{eq-interaction-4B}).
 Thus the proof is complete.$\Box$

\subsection{$L^2$-stability of solutions on the characteristic triangle domain}

Let $\Delta=\Delta(n_1,k_1;k_0)$. We define the following functionals for $(u^{(h)},v^{(h)})$ on $\Delta$.
\begin{definition}\label{def-functional-00}
 For $k_0\le k\le k_1$, define
\[ L_0(k;\Delta)=\sum\limits_{n_1-k_1+k\le n\le n_1+k_1-k} (|u^k_n|^2+|v^k_n|^2),\]
\[D_0(k;\Delta)=\sum_{n_1-k_1+k\le n\le n_1+k_1-k} |u^k_n|^2|v^k_n|^2,\]
\end{definition}
And for the solution $(\widetilde{u}^{(h)},\widetilde{v}^{(h)})$ to (\ref{eq-LBE-B}), we define the followings.
\begin{definition}\label{def-functional-01}
\[
\widetilde{L}_0(k;\Delta)=\sum\limits_{n_1-k_1+k\le n\le n_1+k_1-k} (|\widetilde{u}^k_n|^2+|\widetilde{v}^k_n|^2),\]
\[\widetilde{L}_g(k;\Delta)=\sum\limits_{n_1-k_1+k\le n\le n_1+k_1-k} |g^k_n|^2\]
and
\[ \widetilde{D}_0(k;\Delta)=\sum\limits_{n_1-k_1+k\le n\le n_1+k_1-k} |\widetilde{u}^k_n|^2|\widetilde{v}^k_n|^2.\]
\end{definition}
For the difference $(U^{(h)},V^{(h)})=(\widetilde{u}^{(h)}-u^{(h)},\widetilde{v}^{(h)}-v^{(h)})$, we define the followings.
\begin{definition}\label{def-functional-1}
Let $\Delta=\Delta(n_1,k_1;k_0)$. For $k_0\le k\le k_1$, define
\[
L_1(k;\Delta)=\sum\limits_{n_1-k_1+k\le n\le n_1+k_1-k} (|U^k_n|^2+|V^k_n|^2),\]
\[ Q_1(k;\Delta)=\sum\limits_{n_1-k_1+k\le n\le l\le n_1+k_1-k} \big\{|U^k_n|^2(|v^k_l|^2+|\widetilde{v}^k_l|^2)
+|V^k_l|^2(|u^k_n|^2+|\widetilde{u}^k_n|^2)\big\} \]
and
\[ D_1(k;\Delta)=\sum\limits_{n_1-k_1+k\le n\le n_1+k_1-k} \big\{|U^k_n|^2(|v^k_n|^2+|\widetilde{v}^k_n|^2)
+|V^k_n|^2(|u^k_n|^2+|\widetilde{u}^k_n|^2)\big\}.\]
\end{definition}

\begin{definition}\label{def-functional-2}
For any constant $K>0$ and any $h>0$, define
\[ F_1(k;\Delta)=L_1(k;\Delta)h+KQ_1(k;\Delta)h^2.\]
\end{definition}

To deal with the above functionals, a technical lemma is given as follows.
\begin{lemma}\label{lemma-functional-induction0}
Suppose that  $a^k_n\ge 0$, $b^k_n\ge 0$, $c^k_n\ge 0$ and $d^k_n\ge 0$ for integers $n$ and $k\ge 0$, with
\[ a^{k+1}_{n+1}\le a^k_n+c^k_n \]
and
\[ b^{k+1}_{n+1}\le b^k_n+d^k_n .\]
Given $k_0\ge 0$ and $n_0$, and for $0\le k\le k_0$, let
\[\displaystyle Q_{n_0,k_0}(k)=\sum_{n_0-k_0+k\le n\le l\le n_0+k_0-k}a^k_n b^k_l \]
and
\[\displaystyle D_{n_0,k_0}(k)=\sum_{n_0-k_0+k\le n\le n_0+k_0-k}a^k_n b^k_n.\]
Then for $1\le k\le k_0$ there holds that
\[\displaystyle Q_{n_0,k_0}(k)-Q_{n_0,k_0}(k-1)+D_{n_0,k_0}(k-1)\le E_{n_0,k_0}(k-1),\]
where
\[\displaystyle E_{n_0,k_0}(k)=L_{n_0,k_0}^a(k)L_{n_0,k_0}^d(k)+L_{n_0,k_0}^b(k)L_{n_0,k_0}^c(k)
  \displaystyle +L_{n_0,k_0}^c(k)L_{n_0,k_0}^d(k)\]
and
\[\displaystyle L_{n_0,k_0}^a(k)=\sum_{n_0-k_0+k\le n\le n_0+k_0-k} a^k_n,\]
\[\displaystyle L_{n_0,k_0}^b(k)=\sum_{n_0-k_0+k\le n\le n_0+k_0-k} b^k_n,\]
\[\displaystyle L_{n_0,k_0}^c(k)=\sum_{n_0-k_0+k\le n\le n_0+k_0-k} c^k_n,\]
\[\displaystyle L_{n_0,k_0}^d(k)=\sum_{n_0-k_0+k\le n\le n_0+k_0-k} d^k_n.\]
\end{lemma}
{\it Proof.} For $1\le k\le k_0$, let
\[ Q^{\prime}_{n_0,k_0}(k-1)=\sum_{n_0-k_0+k\le n\le l\le n_0+k_0-k}a^{k-1}_{n-1} b^{k-1}_{l+1}.\]
 Then, using $n^{\prime}=n-1$ and $l^{\prime}=l+1$, we get
 \begin{eqnarray*}
 \displaystyle Q^{\prime}_{n_0,k_0}(k-1) &\displaystyle = \sum^{n_0+k_0-k}_{n=n_0-k_0+k}
 \displaystyle a^{k-1}_{n-1}\big(\sum_{l^{\prime}=n+1}^{n_0+k_0-k+1}b^{k-1}_{l^{\prime}}\big) \\
   &\displaystyle =\sum^{n_0+k_0-k-1}_{n^{\prime}=n_0-k_0+k-1}
  \displaystyle  a^{k-1}_{n^{\prime}}\big(\sum_{l^{\prime}=n^{\prime}+2}^{n_0+k_0-k+1}b^{k-1}_{l^{\prime}}\big) ,
 \end{eqnarray*}
 which gives
 \[ Q^{\prime}_{n_0,k_0}(k-1)+ D_{n_0,k_0}(k-1)\le Q_{n_0,k_0}(k-1).\]
 Therefore,
 \begin{eqnarray*}
 Q_{n_0,k_0}(k)-Q_{ab}(k-1) & \le Q^{\prime}_{n_0,k_0}(k-1)+ E_{n_0,k_0}(k-1) -Q_{n_0,k_0}(k-1)\\
 &\le -D_{n_0,k_0}(k-1)+E_{n_0,k_0}(k-1),
 \end{eqnarray*}
 which completes the proof. $\Box$

 Applying Lemma \ref{lemma-functional-induction0} to the functional $Q_1$ on $\Delta$ yields the following estimates.

\begin{lemma}\label{lemma-functional-induction1}
There exist positive constants $\delta_3$, $h_2$ and $C_4$, independent of  $k_0$, $k_1$ and $n_1$, such that if $h\in (0,h_2)$ and if $L_0(k-1;\Delta)h\le \delta_3$ and $\widetilde{L}_0(k-1;\Delta)h\le \delta_3$, and if $\widetilde{L}_g(k-1,\Delta)h\le \delta_3$, then
\begin{eqnarray}
 Q_1(k;\Delta)-Q_1(k-1;\Delta)+\frac{1}{2}D_1(k-1;\Delta) &\le& C_4 \Big(\widetilde{L}_g(k-1;\Delta)+L_1(k-1;\Delta)\Big) \nonumber\\
   &+& C_4\Lambda(k-1;\Delta) L_1(k-1;\Delta)h. \label{eq-functional-1}
\end{eqnarray}
Here
\[ \Lambda(k-1;\Delta)=D_0(k-1;\Delta)+\widetilde{D}_0(k-1;\Delta)+\widetilde{L}_g(k-1;\Delta) .\]
\end{lemma}
{\it Proof.} In the proof  we fix the domain and omit "$\Delta$" in the functionals, $L_j(k;\Delta)$ and  $D_j(k,\Delta)$, $j=0,1$ etc. for simplification.

For $n_1-k_1+k\le n\le n_1+k_1-k$ and $k_0+1\le k\le k_1$ and for $h\in (0, h_1)$,  Lemma \ref{lemma-interaction-4} gives
\begin{equation}\label{eq-induction-pf1} |U^k_n|^2\le |U^{k-1}_{n-1}|^2+C_3h(|U^{k-1}_{n-1}|^2+|V^{k-1}_{n-1}|^2+E^{k-1}_{n-1} +|g^{k-1}_{n-1}|^2)
\end{equation}
and
\begin{equation}\label{eq-induction-pf2}
 |V^k_l|^2\le |V^{k-1}_{l+1}|^2+C_3h(|U^{k-1}_{l+1}|^2+|V^{k-1}_{l+1}|^2+E^{k-1}_{l+1} +|g^{k-1}_{l+1}|^2)
 \end{equation}
 for $L_0(k-1)h\le \delta_2$ and $\widetilde{L}_0(k-1)h\le \delta_2$ and for $\widetilde{L}_g(k-1,\Delta)h\le \delta_2$;
and  Lemma \ref{lemma-interaction-2} and Lemma \ref{lemma-interaction-2-B} give
\begin{eqnarray} |u^k_n|^2+|\widetilde{u}^k_n|^2 &\le& |u^{k-1}_{n-1}|^2+|\widetilde{u}^{k-1}_{n-1}|^2 +C^{\prime}_2h(|u^{k-1}_{n-1}|^2+|v^{k-1}_{n-1}|^2+|\widetilde{u}^{k-1}_{n-1}|^2+|\widetilde{v}^{k-1}_{n-1}|^2) \nonumber\\
 &\,& +C^{\prime}_2h(|u^{k-1}_{n-1}|^2|v^{k-1}_{n-1}|^2+ |\widetilde{u}^{k-1}_{n-1}|^2|\widetilde{v}^{k-1}_{n-1}|^2+|g^{k-1}_{n-1}|^2) \label{eq-induction-pf3}
\end{eqnarray}
and
\begin{eqnarray} |v^k_l|^2+|\widetilde{v}^k_l|^2 &\le& |v^{k-1}_{l+1}|^2+|\widetilde{v}^{k-1}_{l+1}|^2 +C^{\prime}_2h(|u^{k-1}_{l+1}|^2+|v^{k-1}_{l+1}|^2+|\widetilde{u}^{k-1}_{l+1}|^2+|\widetilde{v}^{k-1}_{l+1}|^2) \nonumber\\
 &\,& +C^{\prime}_2h(|u^{k-1}_{l+1}|^2|v^{k-1}_{l+1}|^2+ |\widetilde{u}^{k-1}_{l+1}|^2|\widetilde{v}^{k-1}_{l+1}|^2+|g^{k-1}_{l+1}|^2) \label{eq-induction-pf4}
\end{eqnarray}
for $L_0(k-1)h\le \delta_2$ and $\widetilde{L}_0(k-1)h\le \delta_2$.
Here $\delta_2$ and  $C^{\prime}_2=C_{\beta}+C_1$ are given as in the proof of Lemma \ref{lemma-interaction-4}.

Then, applying Lemma \ref{lemma-functional-induction0} to the case that $a^k_n=|U^k_n|^2$ and $b^k_n=|v^k_n|^2+|\widetilde{v}^k_n|^2$ and to the case that $a^k_n=|u^k_n|^2+|\widetilde{u}^k_n|^2$ and $b^k_n=|V^k_n|^2$ respectively, we deduce from (\ref{eq-induction-pf1}-\ref{eq-induction-pf4}) that
\begin{eqnarray}
Q_1(k) & \le&  Q_1(k-1)-D_1(k-1)  \nonumber\\
 &+& C_4^{\prime}hL_1(k-1)\big\{ L_0(k-1)+\widetilde{L}_0(k-1)    +
\Lambda(k-1) \big\}
\nonumber\\  &+& C_4^{\prime}\Theta(h,k-1)\{D_1(k-1)
 +\widetilde{L}_g(k-1)\},
  \label{eq-interaction-5}
\end{eqnarray}
where $C^{\prime}_4=(C^{\prime}_2+1)(C_3+1)$ and
\[ \Theta(h,k-1)=\{L_0(k-1)
+\widetilde{L}_0(k-1)\}h+\{D_0(k-1)+\widetilde{D}_0(k-1)+\widetilde{L}_g(k-1)\}h^2.
\]

Now we choose a $\delta_3\in (0, \delta_2)$ so that
\[ -1+ C^{\prime}_4(3\delta_3+2\delta^2_3)\le -\frac{1}{2},\]
and assume that $L_0(k-1)h\le \delta_3$, $\widetilde{L}_0(k-1)h\le \delta_3$ and $\widetilde{L}_g(k-1)h\le \delta_3$,
then
\[ -1+C^{\prime}_4\Theta(h,k-1) \le -\frac{1}{2}.\]
Therefore by (\ref{eq-interaction-5}) and by Lemma \ref{lemma-functional-induction0}, we  get the following,
\begin{eqnarray*}
Q_1(k) & \le&  Q_1(k-1)-\frac{1}{2}D_1(k-1) \\
 &+& C_4^{\prime}L_1(k-1)\big\{ 2\delta    +
\Lambda(k-1)h \big\}
\\  &+&\frac{1}{2} C_4^{\prime}\widetilde{L}_g(k-1),
 \end{eqnarray*}
  which complete the proof.$\Box$

Now we can derive the following estimates on functional $F_1$ by estimates on $D_1$ and $L_1$.
\begin{proposition}\label{prop-functional-1}
 There exist constant $K>0$ and $C_*>0$ independent of $\Delta$ and $h$, such that if $L_0(k_0;\Delta)h\le \delta_*$, $\sup\limits_{k_0\le k\le k_1-1} \widetilde{L}_0(k;\Delta)h\le \delta_*$ and $\sup\limits_{k_0\le k\le k_1-1} \widetilde{L}_g(k;\Delta)h\le \delta_*$, then
\begin{eqnarray}\label{eq-functional-2-0}
 F_1(k;\Delta)-F_1(k-1;\Delta)&\le & C_*h^2\{\widetilde{L}_g(k-1;\Delta)+L_1(k-1;\Delta)\} \nonumber\\
 &+&C_*h^3 \Lambda(k-1;\Delta) L_1(k-1;\Delta)
 \end{eqnarray}
 for $k_0+1\le k\le k_1$ and $h\in (0,h_2)$, where $h_2$ is given by Lemma \ref{lemma-functional-induction1} and
 \[ \Lambda(k-1;\Delta)=D_0(k-1;\Delta)+\widetilde{D}_0(k-1;\Delta)+\widetilde{L}_g(k-1;\Delta).\]
 Moreover,
 \begin{equation}
 F_1(k;\Delta)\le \big\{ F_1(k_0;\Delta)+C_*\sum_{k_0\le j\le k-1} \widetilde{L}_g(j;\Delta)h^2\big\}\exp(\Lambda_*(k_0,k,h))
 \label{eq-functional-2-1}
 \end{equation}
 for $k_0+1\le k\le k_1$ and $h\in (0,h_2)$, where
 \[\Lambda_*(k_0,k,h)=C_*(k-k_0)h+C_*\sum_{j=k_0}^{k-1}(D_0(j;\Delta)+\widetilde{D}_0(j;\Delta)+\widetilde{L}_g(j;\Delta))h^2.\]
\end{proposition}
{\it Proof.} Let $h\in (0, h_2)$. By Lemma \ref{lemma-sum-1}, we have
\begin{equation}\label{eq-functional-2-21} L_0(k;\Delta) \le L_0(k_0;\Delta)\le \delta_*\end{equation} for $k_0+1\le k\le k_1$.

And by Lemma \ref{lemma-interaction-4} we have
\begin{equation}\label{eq-functional-2-2}
L_1(k;\Delta)h\le L_1(k-1;\Delta)h+ 2C_3h^2\big(L_1(k-1;\Delta)+D_1(k-1;\Delta)+\widetilde{L}_g(k-1;\Delta) \big).
\end{equation}
Then we choose a positive constant $K$ so that
\[ \frac{1}{2}K-2C_3>1,\]
which leads to the inequality (\ref{eq-functional-2-0})  by (\ref{eq-functional-2-21})-(\ref{eq-functional-2-2}) and by Lemma \ref{lemma-functional-induction1} directly.

Moreover (\ref{eq-functional-2-0}) implies the following,
\begin{eqnarray}
 F_1(j;\Delta)&\le &
  F_1(j-1;\Delta)\exp(\Lambda_{j-1,j,h})+C_*h^2\widetilde{L}_g(j-1;\Delta) \nonumber\\
 &\le & \big(F_1(j-1;\Delta)+C_*h^2\widetilde{L}_g(j-1;\Delta)\big)\exp(\Lambda_{j-1,j,h}) \label{eq-functional-2-3}
 \end{eqnarray}
 for $k_0+1\le j\le k$,
where
\[ \Lambda_{j-1,j,h}= C_*h+ \big(D_0(j-1;\Delta)+\widetilde{D}_0(j-1;\Delta)+\widetilde{L}_g(j-1;\Delta) \big)h^2.\]
Therefore by induction on $j$, we can  deduce (\ref{eq-functional-2-1}) from (\ref{eq-functional-2-3}). The proof is complete.$\Box$

\section{Compactness of the sequence of solutions in $L^2_{loc}$}

Let $h_2$ be given by Lemma \ref{lemma-functional-induction1}. We consider  the set of solutions $\displaystyle \{(u^{(h)},v^{(h)})\, |\, h\in (0,h_2)\}$ in
 $\displaystyle L^2(R^1\times [0,T])$ for  $T>0$.

\subsection{$L^2$ stability estimates in a strip domain}

Consider the difference
 \[(u^{(h)}(x+\tau,t)-u^{(h)}(x,t), v^{(h)}(x+\tau,t)-v^{(h)}(x,t))\]
 for $|\tau|>0$ and $h\in (0,h_2)$.

 First we have the estimates for such difference in triangle domain $\Delta=\Delta(n_1,k_1;k_0)$ as follows.
\begin{lemma}\label{lemma-stability-0}
Suppose that $\tau=n_0h$ for some integer $n_0\neq 0$ and let
 $\displaystyle (\widetilde{u}^k_n, \widetilde{v}^k_n)=(u^k_{n+n_0}, v^k_{n+n_0})$ for integers $k$ and $n$ with $k\ge 0$. Then,  there exists a constant $C_1(T)>0$, independent of $\tau$ and $\Delta$, such that if $L_0(k_0;\Delta)\le \delta_*$  and if $\widetilde{L}_0(k_0;\Delta)\le \delta_*$ then
 \[ F_1(k; \Delta)\le C_1(T)F_1(k_0;\Delta) \]
 for $k_0\le k\le k_1\le \frac{T}{h}$ and for $h\in (0,h_2)$,
where $\Delta=\Delta(n_1,k_1;k_0)$. Here the constant $\delta_*$ is given by Proposition \ref{prop-functional-1}.
\end{lemma}
{\it Proof.}
  $\displaystyle\{(\widetilde{u}^k_n, \widetilde{v}^k_n)\}_{-\infty<n<\infty, k\ge 0}$ solves the scheme (\ref{eq-LBE-B}) with $g^{k,1}_n=g^{k,2}_n=0$ for $-\infty<n<\infty, k\ge 0$, and $\widetilde{L}_g(k,\Delta)=0$ for $k_0\le k\le k_1$.

 By Lemma \ref{lemma-product-0}, we have
 \[ \sum_{k=k_0}^{k_1}D_0(k; \Delta)h^2\le C_0(T)\]
 and
 \[ \sum_{k=k_0}^{k_1}\widetilde{D}_0(k; \Delta)h^2\le C_0(T)\]
 for some constant $C_0(T)$ depending only on $T$. Here $ k_1\le \frac{T}{h}$.

 Then
  \[\Lambda_*(k_0,k_1,h)\le C_*(k_1-k_0)h+2C_*C_0(T)\le C_*T+2C_*C_0(T),\]
  which, together with proposition \ref{prop-functional-1}, leads to the result. The proof is complete.$\Box$

Lemma \ref{lemma-stability-0} implies the stability of  the solutions in $\Delta$. To extend this result to a strip domain $\{(x,t)|\, x\in R^1, \, 0\le t\le T\}$ for $T>0$, we will divide the strip domain into three  suitable sub-domains and first choose the unbounded domains $\{(x,t)\,|\, |x|\ge A+t, \, 0\le t\le T\}$ for some constant $A>0$ via the following steps.

  \begin{lemma}\label{lemma-stability-2}
  For  $\varepsilon>0$, there exist constants $A_{(\varepsilon)}>0$ and $h_{(\varepsilon)}>0$ such that
  \begin{equation}\label{eq-equiintegrable1}
   \sup_{h\in (0,h_{(\varepsilon)})} \int_{|x|\ge A_{(\varepsilon)}+t} \big( |u^{(h)}(x,t)|^2+|v^{(h)}(x,t)|^2\big)dx\le \varepsilon.
   \end{equation}
  Therefore there exist constants $A>0$ and $h_3\in (0, h_2]$ such that
  \begin{equation}\label{eq-equiintegrable2}
   \sup_{h\in (0,h_3)} \int_{|x|\ge A+t} \big( |u^{(h)}(x,t)|^2+|v^{(h)}(x,t)|^2\big)dx\le \frac{\delta_*}{8}
   \end{equation}
  for  $\delta_*$  given by Proposition \ref{prop-functional-1}. Here $A$ is independent of $h$.
  \end{lemma}
  {\it Proof.} Choose $A_{(\varepsilon)}>0$ so that
\[  \int_{|x|\ge A_{(\varepsilon)}/2} \big( |u_0(x)|^2+|v_0(x)|^2\big)dx\le \frac{\varepsilon}{64}. \]

Then, due to the convergence that  \[ \lim\limits_{h\to 0}(||u^{(h)}(x,0)-u_0||_{L^2(R^1)}+||v^{(h)}(x,0)-v_0||_{L^2(R^1)})=0,\]
 we can choose $h_{(\varepsilon)}>0$ so that
\[ \sup_{h\in (0,h_{(\varepsilon)})} \int_{|x|\ge A_{(\varepsilon)}/2} \big( |u^{(h)}(x,0)|^2+|v^{(h)}(x,0)|^2\big)dx\le \frac{\varepsilon}{4}. \]

Therefore, taking the summation of  (\ref{eq-LBE-conserv-1}) over the domain $\{ (x,s)\,|\, |x|\ge \frac{A_{(\varepsilon)}}{2}+s, \, 0\le s\le t\}$, we have
\[ \int_{|x|\ge A_{(\varepsilon)}+t} \big( |u^{(h)}(x,t)|^2+|v^{(h)}(x,t)|^2\big)dx \le \int_{|x|\ge A_{(\varepsilon)}/2} \big( |u^{(h)}(x,0)|^2+|v^{(h)}(x,0)|^2\big)dx\le \varepsilon \]
for $h\in (0,h_{(\varepsilon)})$, which gives (\ref{eq-equiintegrable1}) and
 (\ref{eq-equiintegrable2}). The proof is complete. $\Box$

And  we can deduce the following stability results in the domain $\{ (x,t)|\, |x|\ge A+2t, \, t>0\}$ by Lemma \ref{lemma-stability-0}.

\begin{lemma}\label{lemma-stability-2-1}
Let $A$ and $h_3$ be the constants given in Lemma \ref{lemma-stability-2} and let $T>0$. Let $n_A=[A/h]+1$. Then  there exists a constant $C(T)>0$ depending on $T>0$ such that
\begin{eqnarray*} \sup_{0\le kh\le T} \int_{|x|\ge 2n_Ah+kh} \big( |u^{(h)}(x+n_0h,kh)-u^{(h)}(x,kh)|^2+|v^{(h)}(x+n_0h,kh)-v^{(h)}(x,kh)|^2\big)dx \\ \le C_2(T)\int_{|x|\ge 2n_Ah} \big( |u^{(h)}(x+n_0h,0)-u^{(h)}(x,0)|^2+|v^{(h)}(x+n_0h,0)-v^{(h)}(x,0)|^2\big)dx
\end{eqnarray*}
for $|n_0|\le \frac{n_A}{2}$ and $h\in (0,h_3)$. Here the constant $C_2(T)$ depends only on $T$.
\end{lemma}
{\it Proof.} Without loss of generality we assume that $A=(n_A-1)h$ and let $\tau=n_0h$. Consider $\Delta(2n_A+2+k_1, k_1,0)$ and $\Delta(-2n_A-2-k_1, k_1,0)$ for $k_1>T/h$.

 It follows from Lemma \ref{lemma-stability-2} that $\int_{2A+ 2h\pm \tau}^{2A+2h+2k_1h\pm\tau}( |u^{(h)}(x,0)|^2+|v^{(h)}(x,0)|^2)dx\le \frac{\delta_*}{8}$ and $\int^{-2A-2h\pm \tau}_{-2A-2h-2k_1h\pm \tau}( |u^{(h)}(x,0)|^2+|v^{(h)}(x,0)|^2)dx\le \frac{\delta_*}{8}$ for $|\tau |\le A/2$.

 Then for $|\tau|\le A/2$ and $h\in (0,h_1)$, by Lemma \ref{lemma-stability-0} we have the following for $0\le k\le T/h$,
 \begin{eqnarray*}  \int_{2A+2h+kh}^{2A+2h+2k_1h-kh} \big( |u^{(h)}(x+\tau,kh)-u^{(h)}(x,kh)|^2+|v^{(h)}(x+\tau,kh)-v^{(h)}(x,kh)|^2\big)dx \\ \le C_2(T)\int_{2A+2h}^{2A+2h+2k_1h} \big( |u^{(h)}(x+\tau,0)-u^{(h)}(x,0)|^2+|v^{(h)}(x+\tau,0)-v^{(h)}(x,0)|^2\big)dx
\end{eqnarray*}
and
\begin{eqnarray*}  \int^{-2A-2h-kh}_{-2A-2k_1h+kh} \big( |u^{(h)}(x+\tau,kh)-u^{(h)}(x,kh)|^2+|v^{(h)}(x+\tau,kh)-v^{(h)}(x,kh)|^2\big)dx \\ \le C_2(T)\int^{-2A-2h}_{-2A-2h-2k_1h} \big( |u^{(h)}(x+\tau,0)-u^{(h)}(x,0)|^2+|v^{(h)}(x+\tau,0)-v^{(h)}(x,0)|^2\big)dx
\end{eqnarray*}
for some constant $C(T)>0$ depending  on $T$,
which lead to the result as $k_1$ goes to infinity. The proof is complete.
$\Box$

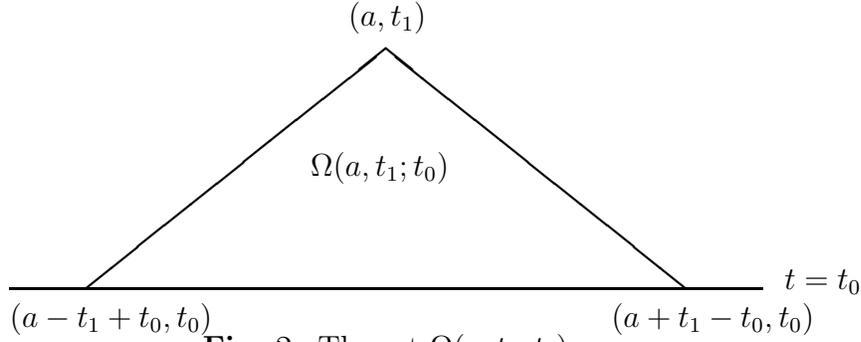
\begin{figure}[h]
\begin{center}
\unitlength=10mm
\begin{picture}(10,3.5)
\thicklines
\put(0,0){\line(1,0){10}}
\put(1,0){\line(5,4){4}}
\put(9,0){\line(-5,4){4}}
\put(0,-0.5){$(a-t_1+t_0,t_0)$}
\put(8,-0.5){$(a+t_1-t_0,t_0)$}\put(10.3,0){$t=t_0$}
\put(4.5,3.5){$(a,t_1)$}
\put(4,1.5){$\Omega(a,t_1; t_0)$}
\end{picture}
\caption{The set $\Omega(a,t_1; t_0)$}\label{fig-domain2}
\end{center}
\end{figure}

Now, for $a\in R^1$ and $t_1\ge t_0\ge 0$, denote
\[ \Omega(a, t_1; t_0)=\{ (x,t)\, |\, a-t_1+t\le x \le a+t_1-t,\, \, t_0\le t\le t_1\}; \]
see Fig. \ref{fig-domain2}.

We consider the $L^2-$stablity estimates in the  domain $\Omega(0,4A+T;0)\cap \Big(R^1\times[0,T]\Big)$. To this end, we first deduce the following for the control of $L^2-$norm of $(h^{(h)},v^{(h)})$ over small intervals.

\begin{lemma}\label{lemma-stability-3}
Let $T>0$ and let $\delta_*$ be the constant given by Proposition \ref{prop-functional-1}. There exist constants $h_4>0$ and $r>0$ such that if $0\le b-a\le 4r$ and $t_0=k_0h\in [0,T]$ for some $k_0>0$  then
\[ \int\limits_{x\in[a,b]} (|u^{(h)}(x,t_0)|^2+|v^{(h)}(x,t_0)|^2)dx \le \frac{\delta_*}{8} \]
for $h\in (0,h_4)$. Here the constants $h_4$  and $r$ are independent of $h$ and $k_0$.
\end{lemma}
{\it Proof.} It suffices to consider two cases: the case that $[a,b]\in (-\infty,-t_0-A)\cup (A+t_0,\infty)$ and the case that
$[a,b]\in (-t_0-4A, 4A+t_0)$. Here $A$ is the constant given by Lemma \ref{lemma-stability-2}.

For case that $[a,b]\in (-\infty,-t_0-A)\cup (A+t_0,\infty)$, the result follows from  Lemma \ref{lemma-stability-2}.

Now we consider the second case that $[a,b]\in (-t_0-4A, 4A+t_0)$.
Due to the hypothesis that \[ \lim\limits_{h\to 0}(||u^{(h)}(x,0)-u_0||_{L^2(R^1)}+||v^{(h)}(x,0)-v_0||_{L^2(R^1)})=0,\]
we can choose a $h_3\in (0, h_1)$ so that
\[ C_1^2( ||u^{(h)}(x,0)-u_0||_{L^2(R^1)}^2+||v^{(h)}(x,0)-v_0||_{L^2(R^1)}^2)\le \delta_*/64 \] for $h\in (0, h_3)$.
In addition we choose a $r\in (0, A/8)$ so that
\[
 C_1^2\int_a^b (|u_0(x)|^2+|v_0(x)|^2) dx\le \delta_*/64
\]
and
\[ C^2_2 T(b-a)\le \delta_*/64\]
for any interval $[a,b]\subset [-T-4A,4A+T]$ with $0<b-a\le 4r$. Here  $h_1$ is the constant given by Lemma \ref{lemma-stability-2}.

Then we have
\begin{equation}\label{eq-small-interval}
\sup_{h\in(0,h_3)}2 C_1^2\int_a^b (|u^{(h)}(x,0)|^2+|v^{(h)}(x,0)|^2) dx+ 4C^2_2 T(b-a)< \delta_*/8
\end{equation}
for any interval $[a,b]\subset [-T-4A,4A+T]$ with $0<b-a\le 4r$.

Noticing that  for $t_0\in [0,T]$ and $[a,b]\subset [t_0-T-4A,T+4A-t_0]$ with $0<b-a\le 4r$, it holds  that $[a-t_0,b-t_0]\subset [-T-4A,T+4A]$ and $[a+t_0,b+t_0]\subset [-T-4A,T+4A]$. Then, by Lemma \ref{lemma-pointwise-1} and by (\ref{eq-small-interval}),
we have
\begin{eqnarray*}
&\,& \int_a^b (|u^{(h)}(x,t_0)|^2+|v^{(h)}(x,t_0)|^2)dx \\ &\le& \int_a^b (C_1|u^{(h)}(x-t_0,0)|+C_2\sqrt{T})^2 dx \\ &\,& +\int_a^b(C_1|v^{(h)}(x+t_0,0)|+C_2\sqrt{T})^2 dx \\
&\le & 2C_1^2\big\{\int_{a+t_0}^{b+t_0} |u^{(h)}(x,0)|^2 dx+\int_{a-t_0}^{b-t_0} |v^{(h)}(x,0)|^2dx\big\} \\ &\,&+4C_2^2T(b-a) < \delta_*/8,
\end{eqnarray*}
which proves the result for second case. The proof is complete.$\Box$

Then, we have the $L^2-$stablity estimates in the  domain $\Omega(0,4A+T;0)\cap \Big(R^1\times[0,T]\Big)$ as follows.

\begin{lemma}\label{lemma-stability-4}
For $T>0$, there exists a constant $C(T)>0$ depending on $T$ only such that
\begin{eqnarray}
\int_{-4A-T+t}^{4A+T-t} \big( |u^{(h)}(x+\tau,t)-u^{(h)}(x,t)|^2+ |v^{(h)}(x+\tau,t)-v^{(h)}(x,t)|^2\big)dx
\nonumber\\
\le C(T) \int^{4A+T}_{-4A-T} \big( |u^{(h)}_0(x+\tau)-u^{(h)}(x)|^2+ |v^{(h)}_0(x+\tau)-v^{(h)}_0(x)|^2\big)dx \label{eq-stability-inner} \,\,
 \, \label{eq-stability-inner}
\end{eqnarray}
for $h\in (0,h_4)$, $|\tau|\le \frac{A}{64}$ and for $t\in [0,T]$. Here $A$ and $h_4$ are  given by Lemma \ref{lemma-stability-2} and  Lemma \ref{lemma-stability-3} respectively.
\end{lemma}
{\it Proof.} Let $r$ and $\delta_*$  be the constants given by Lemma  \ref{lemma-stability-2} and Lemma \ref{lemma-stability-3}. Without loss of generality, we assume that $r=n_rh$ and $T=n_Th$, $A=n_Ah$ for some positive integers $n_r$, $n_T$ and $n_A$, and assume $\displaystyle \frac{T}{4r}=N_0$ and $\displaystyle \frac{A}{r}=N_1$ for some integers $N_0>0$ and $N_1>0$.

Then the proof of the inequality (\ref{eq-stability-inner}) can be carried out by induction on $k$ for $0\le k\le 4N_0$, that is, we assume that (\ref{eq-stability-inner}) holds for $t\in [0,kr]$ and aim to prove that (\ref{eq-stability-inner}) holds for $t\in [0,(k+1)r]$.

Note that
\begin{eqnarray*}
\displaystyle \overline{\Omega(0,4A+T;0)}\cap \{(x,t)\, |\, 0\le t\le (k+1)r\}  \subseteq
\displaystyle \Big(\overline{\Omega(0,4A+T;0)}\cap \{(x,t)\, |\, 0\le t\le kr\} \Big) \cup \\
\displaystyle   \cup_{k+2-4N_0-4N_1\le n\le 4N_0+4N_1-k-2}
\displaystyle  \overline{\Omega(nr, kr+2r ;kr)}.
\end{eqnarray*}

  By Lemma \ref{lemma-stability-0} and Lemma \ref{lemma-stability-3}, there is a constant $C(T)>0$ such that
  \begin{eqnarray}
\int_{(n-k-2)r+t}^{(n+k+2)r-t} \big( |u^{(h)}(x+\tau,t)-u^{(h)}(x,t)|^2+ |v^{(h)}(x+\tau,t)-v^{(h)}(x,t)|^2\big)dx
\nonumber\\
\le C(T) \int^{(n+2)r}_{(n-2)r} \big( |u^{(h)}_0(x+\tau)-u^{(h)}(x)|^2+ |v^{(h)}_0(x+\tau)-v^{(h)}_0(x)|^2\big)dx \nonumber \,\, \, \label{eq-stability-inner-k-n}
\end{eqnarray}
for  $k+2-4N_0-4N_1\le n\le 4N_0+4N_1-k-2$ and $t\in [kr, (k+2)r]$, which leads to the following,
\begin{eqnarray}
\int_{-4A-T+t}^{4A+T-t} \big( |u^{(h)}(x+\tau,t)-u^{(h)}(x,t)|^2+ |v^{(h)}(x+\tau,t)-v^{(h)}(x,t)|^2\big)dx
\nonumber\\
\le C(T,k) \int^{4A+T}_{-4A-T} \big( |u^{(h)}_0(x+\tau)-u^{(h)}(x)|^2+ |v^{(h)}_0(x+\tau)-v^{(h)}_0(x)|^2\big)dx \,\,\, \label{eq-stability-inner-k}
\end{eqnarray}
for  $t\in [kr, (k+2)r]$.
Here the constant $C(T,k)$ depends only on $k$ and $T$ and $r$. Therefore the proof is complete. $\Box$

\subsection{The compactness of the set of the solutions}

To show the compactness of the solution, we consider the difference
 \[(u^{(h)}(\cdot,t),v^{(h)}(\cdot,t))-(u^{(h)}(\cdot-\tau,t-\tau),v^{(h)}(\cdot+\tau,t-\tau))\]
  for $\tau\in R^1$. We have the uniform continuity of $(u^{(h)}(\cdot,t),v^{(h)}(\cdot,t))$ along the characteristic line as follows.
\begin{lemma}\label{lemma-lip-time}
Let $T>0$. For any $\varepsilon>0$, there exist a constant $\delta>0$  such that if  $0\le k_0h\le (k_1+1)h\le T$ with $|k_1-k_0|h<\delta$ and if $h\in (0,\delta)$, then
\[ \sum_{n=-\infty}^{\infty} |u^{k_1+1}_{n+1}-u^{k_0}_{n-k_1+k_0}|^2h < C(T)\big(h+\varepsilon\big)\]
and
\[ \sum_{n=-\infty}^{\infty} |v^{k_1+1}_{n-1}-v^{k_0}_{n+k_1-k_0}|^2h < C(T)\big(h+\varepsilon\big)\] for $h>0$.
Here  $C(T)$ is a constant depending only on $T$.
\end{lemma}
{\it Proof.}
 For $0\le j\le k_1-k_0+1$ and $-\infty<n<\infty$, the first equation in (\ref{eq-LBE}) gives the following
\begin{eqnarray*}
(1-\frac{ih\alpha}{2} |v^{k_1-j}_{n-j}|^2)(u^{k_1+1-j}_{n+1-j}-u^{k_0-j}_{n-j}) &=&\frac{imh}{2}(v^{k_1+1-j}_{n-1-j}+v^{k_1-j}_{n-j}) +ih\alpha u^{k_1-j}_{n-j}|v^{k_1-j}_{n-j}|^2 \\ &\,& \displaystyle +\frac{ih\beta}{2}(v^{k_1+1-j}_{n-1-j}+v^{k_1-j}_{n-j})G(u^{k_1-j}_{n-j},v^{k_1-j}_{n-j}),
\end{eqnarray*}
 which leads to
\begin{eqnarray}
|u^{k_1+1-j}_{n+1-j}-u^{k_1-j}_{n-j}|&\le&\frac{mh}{2}|v^{k_1+1-j}_{n-1-j}+v^{k_1-j}_{n-j}| +h|\alpha | |u^{k_1-j}_{n-j}||v^{k_1-j}_{n-j}|^2 \nonumber\\ &\,& +h|\beta|(|v^{k_1+1-j}_{n-1-j}|+|v^{k_1-j}_{n-j}|)|u^{k_1-j}_{n-j}||v^{k_1-j}_{n-j}|. \label{eq-LBE-j}
\end{eqnarray}

Now taking the summation of (\ref{eq-LBE-j}) over $j$ for $0\le j\le k_1-k_0+1$, we have
\begin{eqnarray*}
|u^{k_1+1}_{n+1}-u^{k_0}_{n_0}|&\le& \sum\limits_{0\le j\le k_1-k_0+1}\frac{mh}{2}(|v^{k_1+1-j}_{n-1-j}|+|v^{k_1-j}_{n-j}|) \\ &\,&  +\sum\limits_{0\le j\le k_1-k_0+1} h(|\alpha|+|\beta|)(|v^{k_1+1-j}_{n-1-j}|+|v^{k_1-j}_{n-j}|)|u^{k_1-j}_{n-j}||v^{k_1-j}_{n-j}| ,
\end{eqnarray*}
where $n_0=n_1-k_1+k_0$.
Then,
\begin{eqnarray}
|u^{k_1+1}_{n+1}-u^{k_0}_{n_0}|^2&\le& m^2h^2\big\{\sum\limits_{0\le j\le k_1-k_0+1}(|v^{k_1+1-j}_{n-1-j}|+|v^{k_1-j}_{n-j}|)\big\}^2 \nonumber\\ &\,& + 2h^2(|\alpha|+|\beta|)^2\big\{\sum\limits_{0\le j\le k_1-k_0+1}(|v^{k_1+1-j}_{n-1-j}|+|v^{k_1-j}_{n-j}|)|u^{k_1-j}_{n-j}||v^{k_1-j}_{n-j}| \big\}^2 \nonumber\\
&\le& m^2h^2\sum_{j=0}^{k_1-k_0+1}(|v^{k_1+1-j}_{n-1-j}|^2+|v^{k_1-j}_{n-j}|^2)(k_1-k_0+1) \nonumber \\
&+& C_{\alpha,\beta}h^2\sum_{j=0}^{k_1-k_0+1}(|v^{k_1+1-j}_{n-1-j}|^2 +|v^{k_1-j}_{n-j}|^2)\sum_{j=0}^{k_1-k_0+1}|u^{k_1-j}_{n-j}||v^{k_1-j}_{n-j}|^2 \label{eq-LBE-jSum0}
\end{eqnarray}
for $C_{\alpha,\beta}=8(|\alpha|+|\beta|)^2$.

Next we deal with three terms in the last inequality in (\ref{eq-LBE-jSum0}). First  Corollary \ref{coro-conservation} gives
\begin{eqnarray}\label{eq-LBE-jSum1}
\sum_{n=-\infty}^{\infty}\sum_{j=0}^{k_1-k_0+1}(|v^{k_1+1-j}_{n-1-j}|^2+|v^{k_1-j}_{n-j}|^2)h^2 &\le& \sum_{j=0}^{k_1-k_0+1} \sum_{n=-\infty}^{\infty} (|v^{k_1+1-j}_{n-1-j}|^2+|v^{k_1-j}_{n-j}|^2) h^2\nonumber\\
 &=& \sum_{j=0}^{k_1-k_0+1} \sum_{n=-\infty}^{\infty} (|v^{k_1+1-j}_{n}|^2+|v^{k_1-j}_{n}|^2)h^2 \nonumber\\
&\le& 2(k_1-k_0+1)h\sum_{n=-\infty}^{\infty} (|u^{0}_{n}|^2+|v^{0}_{n}|^2)h,\,\,\qquad\,\,
\end{eqnarray}
and Lemma \ref{lemma-product-0} gives
\begin{equation}\label{eq-LBE-jsum3}
\sum^{+\infty}_{n=-\infty}|u^{k_1-j}_{n-j}|^2|v^{k_1-j}_{n-j}|^2h^2\le C_0^2+(4C_1^2+C_1)C_0T^2.
\end{equation}

Moreover  due to the convergence in $L^2(R^1)$ of the sequence $\{(u^{(h)}(x,0),v^{(h)}(x,0))\}$, there is a constant $\delta>0$  such that
\begin{equation}
\int^{(n+k_1-k_0+2)h}_{(n-k_1+k_0-2)} (|u^{(h)}(x-k_0h,0)|^2+|v^{(h)}(x+k_0h,0)|^2)dx \le \varepsilon
\end{equation}
for $(k_1-k_0)h<\delta$ and $h\in (0,\delta)$.

Then for last terms in (\ref{eq-LBE-jSum0}),
Lemma \ref{lemma-sum-1} and  Lemma \ref{lemma-pointwise-1} gives
\begin{eqnarray}\label{eq-LBE-jSum2}
&\,&\sum_{j=0}^{k_1-k_0+1}(|v^{k_1+1-j}_{n-1-j}|^2+|v^{k_1-j}_{n-j}|^2)h  \nonumber\\&\le& \sum_{l=n-k_1+k_0-2}^{n+k_1-k_0+2} (|u^{k_0}_l|^2+|v^{k_0}_l|^2)h \nonumber\\
&\le& \sum_{l=n-k_1+k_0-2}^{n+k_1-k_0+2} \{2C_1 (|u^0_{l-k_0}|^2+|v^0_{l+k_0}|^2)h+4C_4k_0h^2\} \nonumber \\
&\le & 2C_1\int^{(n+k_1-k_0+2)h}_{(n-k_1+k_0-2)} (|u^{(h)}(x-k_0h,0)|^2+|v^{(h)}(x+k_0h,0)|^2)dx \nonumber \\
&\,& +8C_1C_4Th \nonumber \\
&\le& 2C_1\varepsilon +8C_1C_4Th,
\end{eqnarray}
Therefore  the result can be deduced from (\ref{eq-LBE-jSum0}) by (\ref{eq-LBE-jSum1}), (\ref{eq-LBE-jsum3}) and (\ref{eq-LBE-jSum2}). The proof is complete. $\Box$

The above lemma has a equivalent one as follows.
\begin{lemma}\label{lemma-stability-space}
Let $T>0$. For any $\varepsilon>0$, there exist a constant $\delta>0$ such that if $0<t_0<t_1$ with $|t_0-t_1|<\delta$ and if $h\in (0, \delta)$ then
\[
\int^{\infty}_{-\infty}|u^{(h)}(x,t_1)-u^{(h)}(x-t_1+t_0, t_0)|^2\le C(T)(h+\varepsilon)
\] and
\[ \int^{\infty}_{-\infty}|v^{(h)}(x,t_1)-v^{(h)}(x+t_1-t_0, t_0)|^2\le C(T)(h+\varepsilon)
\]
for $h>0$. Here the constant $C(T)$  depends only on $T$.
\end{lemma}

As a consequence of Lemma \ref{lemma-stability-2}, Lemma \ref{lemma-stability-3}, Lemma \ref{lemma-stability-4} and Lemma \ref{lemma-stability-space}, we can get directly the compactness property of $\{(u^{(h)}, v^{(h)})\}$ as follows.

\begin{proposition}\label{prop-compctness}
Let $(u_0,v_0)\in L^2(R^1)$. Then, for any sequence $\{h_l\}_{l=0}^{\infty}$ with $h_l>0$ for $l\ge 0$ and $\displaystyle \lim_{l\to\infty} h_l=0$, the sequence $\{ (u^{(h_l)},v^{(h_l)}) \}$ is relatively compact in $C([0,T]; L^2(R^1))$ for any $T>0$.
\end{proposition}

\section{Uniqueness of limit and proof of the main result}

Our aim is to show that the sequence $\{ u^{(h)}, v^{(h)}\}$ is strongly convergent in $L^2$ to the unique solution to the problem (\ref{eq-dirac}) and (\ref{eq-dirac-initialv}) as $h$ goes to zero. To this end, we first recall the result in \cite{zhang-zhao} on the well-posedness of global strong solution to (\ref{eq-dirac}) and (\ref{eq-dirac-initialv}).
\begin{theorem}\label{thm-global}
For $(u_0(x), v_0(x))\in L^2(R^1)$, the Cauchy problem (\ref{eq-dirac}) and (\ref{eq-dirac-initialv}) has a unique global strong solution $(u_*,v_*)\in C([0,\infty); L^2(R^1))$ . Moreover, $|u_*||v_*|\in L^2(R^1\times[0,T])$ for any $T>0$.
\end{theorem}

More precisely, according to \cite{zhang-zhao},  there exists a sequence of smooth solution $(u_{j},v_{j})$ to (\ref{eq-dirac}) with $(u_{j}(x,0), v_{j}(x,0))\in C^{\infty}_c(R^1)$ such that
\begin{equation}
\lim_{j\to\infty} \max_{0\le t\le T}\big(||u_*(\cdot,t)-u_j(\cdot,t)||_{L^2(R^1)}+||v_*(\cdot,t)-v_j(\cdot,t)||_{L^2(R^1)}\big)= 0
\end{equation}
 and
\begin{equation}
\lim_{j\to\infty} \big(||u_*v_*-u_jv_j||_{L^2(R^1\times [0,T])}+||u_*\overline{v_*}-u_j\overline{v_j}||_{L^2(R^1\times [0,T])}\big)= 0
\end{equation}
for any $T>0$.

And by the convergence of the sequence $\{(u_j,v_j)\}$, we can assume that
\[ \sup_j \int_{|x|\ge A} (|u_j(x,0)|^2+|v_j(x,0)|^2)dx\le \frac{\delta_*}{7}, \]
where the constant $A$ is given by Lemma \ref{lemma-stability-3}. Then, it is proved in \cite{zhang-zhao} by the conservation of the charge that the followings hold.
\begin{lemma}\label{lemma-sequence-smooth-1}
 For $t\ge 0$, there holds that
\[ \sup_j \int_{|x|\ge A+t} (|u_j(x,t)|^2+|v_j(x,t)|^2)dx\le \frac{\delta_*}{4}.\]
\end{lemma}

\begin{lemma}\label{lemma-sequence-smooth-2}
Let $T>0$. There exists a $r^{\prime}>0$ such that if $t-(A+2T)\le a\le b\le (A+2T)-t$  with $|a-b|\le 16r^{\prime}$ and $t\in[0,T]$ then
\[ \sup_j \int^b_a (|u_j(x,t)|^2+|v_j(x,t)|^2)dx\le \frac{\delta_*}{4}.\]
\end{lemma}

In the next we assume that $r^{\prime}=r$ for simplification, and consider the difference between the QLB solutions $(u^{(h)},v^{(h)})$ and the smooth solution $(u_j,v_j)$.

 Let
\[ (u_j^{(h)},v_j^{(h)})(x,t)=(u_j,v_j)(nh,kh), \quad nh\le x<(n+1)h, \, \, kh\le t <(k+1)h \]
for $k\ge 0$ and $-\infty<n<\infty$, and denote
\[ (u_{j,n}^k, v_{j,n}^k)=(u_j,v_j)(nh,kh).\]
Then
\begin{eqnarray}\label{eq-LBE-B-ju}
u^{k+1}_{j,n+1}-u^k_{j,n}=\frac{imh}{2}(v^{k+1}_{j,n-1}+v^k_{j,n}) +\frac{i\alpha h(u^{k+1}_{j,n+1}+u^k_{j,n})}{2}|v^k_{j,n}|^2  \nonumber
 \\  +\frac{ih\beta}{2}(v^{k+1}_{j,n-1}+v^k_{j,n})G(u^k_{j,n},v^k_{j,n}) +g^{k,1}_{j,n}h
\end{eqnarray}
and
\begin{eqnarray}\label{eq-LBE-B-jv}
 v^{k+1}_{j,n-1}-v^k_{j,n} =\frac{imh}{2}(u^{k+1}_{j,n+1}+u^k_{j,n}) +\frac{i\alpha h(v^{k+1}_{j,n-1}+v^k_{j,n})}{2}|u^k_{j,n}|^2 \nonumber
  \\ +\frac{ih\beta}{2}(u^{k+1}_{j,n+1}+u^k_{j,n})G(u^k_{j,n},v^k_{j,n}) +g^{k,2}_{j,n}h,
\end{eqnarray}
where
\begin{eqnarray*} g^{k,1}_{j,n}=\int^1_0 \big( im v_j+iN_1(u_j,v_j)\big)((n+\tau)h, (k+\tau)h) d\tau
-\frac{im}{2}(v^{k+1}_{j,n-1}+v^k_{j,n})
 \\ -\frac{i\alpha (u^{k+1}_{j,n+1}+u^k_{j,n})}{2}|v^k_{j,n}|^2
-\frac{i\beta}{2}(v^{k+1}_{j,n-1}+v^k_{j,n})G(u^k_{j,n},v^k_{j,n}) \end{eqnarray*}
and
\begin{eqnarray*} g^{k,2}_{j,n}=\int^1_0 \big( im u_j+iN_2(u_j,v_j)\big)((n+\tau)h, (k+\tau)h) d\tau
-\frac{im}{2}(u^{k+1}_{j,n+1}+u^k_{j,n})
\\ -\frac{i\alpha (v^{k+1}_{j,n-1}+v^k_{j,n})}{2}|u^k_{j,n}|^2 \nonumber
  -\frac{i\beta}{2}(u^{k+1}_{j,n+1}+u^k_{j,n})G(u^k_{j,n},v^k_{j,n}).
\end{eqnarray*}

Direct computation gives the following.
\begin{lemma}\label{lemma-sequence-smooth-j}
There holds that
\[ \max_{k\ge0, -\infty<n<\infty} |g^{k,l}_{j,n}|\le M_j(T)h \]
for $h>0$, $l=1,2$ and $j=1,2,\cdots$, where $M_j(T)=C(m,\alpha,\beta)M_j^1(T)M_j^0(T)^2$ and $C(m,\alpha,\beta)=6(m+|\alpha|+|\beta|)$,
\[ M_j^1(T)=\big(\max_{R^1\times [0,T]} |u_{jt}|+\max_{R^1\times[0,T]}|u_{jx}|+ \max_{R^1\times [0,T]} |v_{jt}|+\max_{R^1\times[0,T]}|v_{jx}|\big)\]
and
\[M_j^0=\big(1+\max_{R^1\times[0,T]}|u_{j}|+ \max_{R^1\times [0,T]} |v_{j}|\big).\]
Moreover, for $j\ge 0$, there exists a $h_{*,j}>0$ such that  if $0\le kh\le T$ then
\begin{equation}
\sup_{h\in (0, h_{*,j})}\sum_{n=-\infty}^{\infty} (|g^{k,1}_{j,n}|^2+|g^{k,2}_{j,n}|^2)\le \frac{\delta_*}{4}.
\end{equation}
and
\begin{equation}
\sup_{h\in (0,h_{*,j})}\sum_{n=-\infty}^{\infty} \sum_{p=0}^{k_T}(|g^{p,1}_{j,n}|^2+|g^{p,2}_{j,n}|^2)<\infty.
\end{equation}
Here $k_T=[T/h]+1$.
\end{lemma}

{\bf Proof of Theorem \ref{thm-main}.} Due to the Proposition \ref{prop-compctness}, it remains to prove that the  strong solution $(u_*,v_*)$ of  the problem (\ref{eq-dirac}) and (\ref{eq-dirac-initialv}) is the limit of any convergent subsequence of $\{(u^{(h)},v^{(h)})\}$  as $h$ goes to zero.

Let $\{ (u^{(h_l)},v^{(h_l)})\}$ be a subsequence of $\{(u^{(h)},v^{(h)})\}$ with
$\displaystyle\lim_{l\to \infty}h_l=0$. Then by Proposition \ref{prop-compctness}, it has a subsequence which is convergent in $L^2(R^1\times [0,T])$ to a $(u^{\prime},v^{\prime})$ for any $T>0$. We  still denoted this convergent subsequence by $\{ (u^{(h_l)},v^{(h_l)})\}$ for simplification.

Our aim is to show that $(u_*,v_*)=(u^{\prime},v^{\prime})$. To this end, we divide the time interval $[0,T]$ by the points $t=pr$, $p=0,1,\cdots, k_T$. Here we assume that $K_Tr=T$ for some integer $k_T\ge 0$.

 Now we use the induction on $p$, that is, we assume that $(u_*,v_*)(x,t)=(u^{\prime},v^{\prime})(x,t)$ for $(x,t)\in R^1\times [0,pr]$ for $0\le p\le k_T$.It suffices to consider the case that $p<k_T$.

By Lemma \ref{lemma-sequence-smooth-1}, Lemma \ref{lemma-sequence-smooth-2} and Lemma \ref{lemma-sequence-smooth-j}, we applied Proposition \ref{prop-functional-1} to $(u^{(h)},v^{(h)})$ and $(\widetilde{}u^{(h)}, \widetilde{v}^{(h)})=(u^{(h)}_j,v^{(h)}_j)$ to get the following on $\Omega(a, (p+2)r; pr)$ for any $a\in R^1$,
\begin{eqnarray}\label{eq-sequence-smooth-j1}
\sup_{pr\le t\le (p+2)r} \int_{a-(p+2)r+t}^{a+(p+2)r-t} \big(|(u^{(h)}-u^{(h)}_j)(x,t)|^2+ |(v^{(h)}-v^{(h)}_j)(x,t)|^2\big)dx \nonumber\\
\le C(T) \int_{a-2r}^{a+2r} \big(|(u^{(h)}-u^{(h)}_j)(x,pr)|^2+ |(v^{(h)}-v^{(h)}_j)(x,pr)|^2\big)dx \nonumber \\
+C(T)\Big(\sum_{k=0}^{k_T}\sum^{\infty}_{n=-\infty}(|g^{k,1}_{j,n}|^2+|g^{k,2}_{j,n}|^2)h^2\Big)\exp\big(C_*(T)+C_*
\sum_{k=0}^{k_T}\sum^{\infty}_{n=-\infty}(|g^{k,1}_{j,n}|^2+|g^{k,2}_{j,n}|^2)h\big).
\end{eqnarray}
Since $(u_j,v_j)\big|_{R^1\times[0,T]}$ has compact support in $R^1\times [0,T]$, then by Lemma \ref{lemma-sequence-smooth-j}, we put $h=h_l$ in (\ref{eq-sequence-smooth-j1}) and take the limit as $h_l$ goes to $0$ to deduce that
\begin{eqnarray}\label{eq-sequence-smooth-j2}
\sup_{pr\le t\le (p+2)r} \int_{a-(p+2)r+t}^{a+(p+2)r-t} \big(|(u^{\prime}-u_j)(x,t)|^2+ |(v^{\prime}-v_j)(x,t)|^2\big)dx \nonumber\\
\le C(T) \int_{a-2r}^{a+2r} \big(|(u^{\prime}-u_j)(x,pr)|^2+ |(v^{\prime}-v_j)(x,pr)|^2\big)dx
\end{eqnarray}
for $j\ge 1$.

Then, we can take the limit (\ref{eq-sequence-smooth-j2}) as $j$ goes to $\infty$ to conclude that $(u_*,v_*)(x,t)=(u^{\prime},v^{\prime})(x,t)$ on $\Omega(a, (p+2)r; pr)$ for any $a\in R^1$. Therefore
$(u_*,v_*)(x,t)=(u^{\prime},v^{\prime})(x,t)$ on $R^1\times [0, (p+1)r]$.

Thus carrying out the induction steps yields that $(u_*,v_*)(x,t)=(u^{\prime},v^{\prime})(x,t)$ on $R^1\times [0, T]$ for any $T>0$. The proof is complete.$\Box$


\end{document}